\documentclass[10pt,a4paper]{article}
\usepackage[T1]{fontenc}
\usepackage[utf8]{inputenc}
\usepackage[english]{babel}
\usepackage{microtype} 
\usepackage{xcolor}
\usepackage{graphicx} 
\usepackage{caption}
\usepackage{subfig}
\usepackage{overpic}
\usepackage{geometry}

\usepackage{amsmath}
\usepackage{amssymb}
\usepackage{amsfonts}
\usepackage{mathtools}
\usepackage{amsthm} 
\usepackage{mathrsfs} 
\usepackage{hyperref}
\hypersetup{
    colorlinks=true,
    linkcolor=blue,
    linktocpage=true,
    citecolor=purple,
    anchorcolor=black,
    filecolor=magenta,      
    urlcolor=blue,
    pdfstartview=FitH,
    pdftitle={Octahedrally and Icosahedrally Symmetric Solutions of the N-Vortex Problem on the Sphere},
    pdfauthor={Pietro Pomaro},
    pdfkeywords={point vortices, $N$-vortex problem on the sphere, symmetric solutions, equilibrium and collision configurations, invariant theory of reflection groups}
}

\theoremstyle{definition}

\newtheorem{defn}{Definition}[section]
\newtheorem*{defn*}{Definition}

\theoremstyle{plain}
\newtheorem{thm}[defn]{Theorem}
\newtheorem*{thm*}{Theorem}

\newtheorem*{lemma*}{Lemma}

\newtheorem{prop}[defn]{Proposition}
\newtheorem*{prop*}{Proposition}

\newtheorem*{cor*}{Corollary}

\theoremstyle{remark}

\newtheorem{rem}[defn]{\textsc{Remark}}
\newtheorem*{rem*}{\textsc{Remark}}

\newtheorem*{esempio*}{Example}

\newtheoremstyle{colonstyle}
  {3pt}
  {3pt}
  {}
  {}
  {\bfseries}
  {:}
  {0.5em}
  {}
\theoremstyle{colonstyle}

\newtheorem{auxnamedthm}{\nameofthethm}%
{\end{auxnamedthm}}

\newtheorem*{auxnamedthm*}{\nameofthethm}%
\newenvironment{namedthm*}[1]{%
\newcommand{\nameofthethm}{#1}\begin{auxnamedthm*}}%
{\end{auxnamedthm*}}

\usepackage{fancyhdr}
\fancypagestyle{plain}{
\fancyhead{}
\fancyfoot{}
\fancyhead[LO]{}
\fancyhead[CO]{}
\fancyhead[RO]{}
\fancyfoot[LO]{}
\fancyfoot[CO]{}
\fancyfoot[RO]{\thepage}
\fancyhead[LE]{}
\fancyhead[CE]{}
\fancyhead[RE]{}
\fancyfoot[LE]{\thepage}
\fancyfoot[CE]{}
\fancyfoot[RE]{}

}

\usepackage[autostyle, english=british]{csquotes}
\usepackage[backend=biber, style=numeric-comp, sorting=nyt]{biblatex}
\newcommand\restr[2]{{
  \left.\kern-\nulldelimiterspace 
  #1 
  \vphantom{\big|}
  \right|_{#2}
  }}
  
\newcommand{\numberset}{\mathbb}

\newcommand{\R}{\numberset{R}}

\renewcommand{\phi}{\varphi}
\renewcommand{\theta}{\vartheta}
\renewcommand{\epsilon}{\varepsilon}

\DeclarePairedDelimiter{\abs}{\lvert}{\rvert}
\DeclarePairedDelimiter{\norm}{\lVert}{\rVert}

\title{Octahedrally and Icosahedrally Symmetric Solutions of the $N$-Vortex Problem on the Sphere}
\author{Pietro Pomaro\thanks{Dipartimento di Matematica “Tullio Levi-Civita”, Università di Padova, Via Trieste 63, 35121 Padova, Italy. ptrpmr01@gmail.com}}
\date{October 2026}

\begin{document}

\maketitle
\pagestyle{plain}
\begin{abstract}
We investigate the octahedrally and icosahedrally symmetric solutions of the $N$-vortex problem on the sphere, for $N = 24$ and $N = 60$ respectively, in the case of equal vorticities. Exploiting the symmetries arising from rotations and permutations of the vortices, the problem is reduced to the analysis of a one-degree-of-freedom Hamiltonian system on the sphere. The original contributions of this article concern the analysis of this reduced system yielding in particular a classification of the equilibrium and collision configurations with the prescribed symmetries. The computational complexity of the analysis is overcome by exploiting the associated discrete symmetries through the invariant theory of reflection groups and Palais' principle of symmetric criticality, and, in some cases, by resorting to numerical approximations.
\end{abstract}

    
\tableofcontents
\section{Introduction}

The problem of $N$-point vortices moving on a sphere is by now a classical topic in fluid dynamics. The foundations were laid by the work of Gromeka \cite{Gro} in 1885 and by Zermelo's dissertation \cite{Zer} in 1902. Later, in 1977, Bogomolov \cite{Bog} derived the equations of motion. A comprehensive introduction to the subject (reflecting developments up to the 2000s) is presented in Newton's book \cite{New} and in the references therein. A non-exhaustive list of modern references includes \cite{PolDri, KidNew, LimMonRob, BorKilMam, LauMonRob, CalDoeGar, ModViv, Ohs, CGGL}. 

In this article, we study families of octahedrally and icosahedrally symmetric solutions of the $N$-vortex problem on the sphere, in the case of equal vorticities, for $N=24$ and $N=60$. The analysis of these families was proposed in the “Future Work” section of \cite{AzNa}.
These two values of $N$ are dictated by the symmetry: a symmetric configuration is the orbit of a single vortex under a finite rotation group $K<SO(3)$, hence consists of $\abs{K}$ vortices, and we take $K$ to be the chiral octahedral group $\mathbb{O}$ ($\abs{\mathbb{O}}=24$) and the chiral icosahedral group $\mathbb{I}$ ($\abs{\mathbb{I}}=60$), identified with appropriate representations in $SO(3)$. Such solutions are characterised by the property that the motion of the entire system can be recovered from the motion of a single vortex through the action of $\mathbb{O}$ and $\mathbb{I}$, respectively. In particular, our analysis focuses on the classification of equilibrium and collision configurations with these symmetries.

\subsection{Original contributions}

Combining algebraic techniques and a numerical investigation, we are able to obtain the following original results:
\begin{itemize}
\item{Classification of octahedrally symmetric equilibrium and collision configurations (see Theorem \ref{thm: Oct eqs})};
\item{Classification of icosahedrally symmetric collision configurations and  existence of icosahedrally symmetric equilibrium configurations (see Theorem \ref{thm: Ico eqs})}.
\end{itemize} On the basis of our analysis, we conjecture that the configurations in item (i) of Theorem \ref{thm: Ico eqs} include all possible icosahedrally symmetric equilibrium configurations, providing a complete classification.
Within the framework of \cite{AzNa}, the classification of the symmetric equilibrium and collision configurations reduces to the classification of the critical points of the restriction to the sphere $\mathbb{S}^2$ of appropriate polynomials $p_{\mathbb{O}}$ and $p_{\mathbb{I}}$ in $\R[x,y,z]$. However, this classification is very challenging even with numerical methods, as these polynomials have degrees $46$ for $N=24$ and $118$ for $N=60$. We address this difficulty by taking advantage of the invariance of $p_{\mathbb{O}}$ and $p_{\mathbb{I}}$ under the action of the full octahedral group $\mathbb{O}_h$ and the full icosahedral group $\mathbb{I}_h$. In both cases, such an invariance allows one to restrict the analysis to a fundamental domain $\mathcal{R}\subset\mathbb{S}^2$ for the action. The interior $\mathrm{int}(\mathcal{R})$ can be parametrised through the values of suitable invariant polynomials (the fundamental invariants, whose existence is guaranteed by the Chevalley--Shephard--Todd theorem), and this choice of parametrisation considerably lowers the degree of $p_{\mathbb{O}}$ and $p_{\mathbb{I}}$, easing the search for critical points. On the other hand, Palais' principle of symmetric criticality is applied to the analysis of the boundary $\partial\mathcal{R}$.

\subsection{Structure of the article}

The second section is devoted to recalling the framework of \cite{AzNa} to analyse solutions of the $N$-vortex problem on the sphere with prescribed symmetries. The discussion is not carried out in its full generality, since it is tailored to our specific goals in the following sections.
In the third section, we explain the strategy to determine the critical points of the (regularised) reduced Hamiltonian, discussing our choice of the fundamental regions and their parametrisation through the invariants. The strategy is then applied in the fourth and fifth sections to establish our classification results in the octahedral and icosahedral cases, respectively. We present detailed descriptions and illustrations of the symmetric equilibrium and collision configurations and of the (reduced) phase space of symmetric solutions. Appendix A contains a description of the polyhedra corresponding to equilibrium configurations, whereas Appendix B contains the complete expressions of polynomials and objects previously introduced that are too long or inconvenient to be included in the body of the article.

\begin{namedthm*}{\textbf{Note}} The computer algebra system \textit{Mathematica} was used to draw all the figures and to help factorise large polynomials and provide the numerical values of their roots. Notebooks can be provided upon reasonable request.\end{namedthm*}

\section{Preliminaries}

This section contains the necessary material to begin our analysis. We present a brief discussion and refer to \cite{AzNa} for details and proofs.

\subsection{The equations of motion and their Hamiltonian structure}

Consider the manifold $M \coloneqq \mathbb{S}^2\times\dots\times\mathbb{S}^2\subset (\R^3)^N$ given by the product of $N$ copies of $\mathbb{S}^2\coloneqq \{r\in\R^3\,\colon\,\norm{r}=1\}
$, endowed with the symplectic form \begin{equation}\label{form}\Omega=\sum_{i=1}^N \pi_i^\ast\omega_{\,\mathbb{S}^2}\,,\end{equation} where $\left(\omega_{\,\mathbb{S}^2}\right)_r(u,v)=r\cdot(u\times v)$, $u,v\in T_r\mathbb{S}^2$ is the standard area 2-form on $\mathbb{S}^2$ and $\pi_i$ is the projection onto the $i$-th factor.
Under the assumption of equal vorticity, after a suitable time rescaling, the motion of the $N$ vortices on $\mathbb{S}^2$ is described by the Hamiltonian:
\begin{equation}\label{Ham}
    H(\mathbf{r})=-\frac{1}{2}\sum_{1\le i<j\le N}\ln \big(\norm{r_j-r_i}^2\big),\quad \mathbf{r}=(r_1,\dots,r_N)\in M\,.
\end{equation}The Hamiltonian vector field $X_H$, defined by $dH=i_{X_H}\Omega$, yields the following equations of motion\footnote{Here $\nabla H(\mathbf{r})$ denotes the gradient (with respect to the standard Euclidean metric $\big(\,(\R^{3})^N,\,\cdot\,\big)$) of any extension of $H$ to an open set $U\subset(\R^{3})^N$ containing $M$.}:
\begin{equation} \dot{r}_j=-r_j\times(\nabla H(\mathbf{r}))_j=\sum_{\substack{{i=1}\\{i\ne j}}}^N \frac{r_i\times r_j}{\norm{r_i-r_j}^2}\,, \quad j=1,\dots,N \,.\label{system}
\end{equation} Both the Hamiltonian $H$ and the equations of motion are undefined at the \textit{collision set}
\[\mathscr{C}= \left\{\,(r_1,...,r_N) \in M : r_i = r_ j\,\text{ for some } \,i \ne j \,\right\}\,,\] consisting of configurations where at least two vortices occupy the same position. It can be proven that the flow of \eqref{system} on $M\setminus\mathscr{C}$ is complete, which means that any configuration without collisions cannot evolve into one with collisions.

\subsection{Symmetries and symmetric solutions}

The Hamiltonian \eqref{Ham} is invariant under the symplectic diagonal action of $SO(3)$ rotating all the vortices. The main reason for assuming that all vortices have the same vorticity is that it ensures an additional (discrete) symmetry: the permutation group $S_N$ acts on $M$ by permuting the vortices.
It is possible to combine the two actions into a single symplectic action of the direct product $\hat{G}\coloneqq S_N\times SO(3)$. The group $\hat{G}$ acts on $M$ by: 
\[\Psi^{\hat{G}}\colon \hat{G}\times M\to M,\quad ((\sigma,\,g),\,\mathbf{r})\mapsto \Psi^{\hat{G}}_{(\sigma,\,g)}(\mathbf{r})= ( g\,r_{\sigma^{-1}(1)},\dots,g\,r_{\sigma^{-1}(N)})\,.\]
This is a symplectic action and the Hamiltonian $H$ of \eqref{Ham} is invariant with respect to it. In particular, the equations of motion \eqref{system} are $\hat{G}$-equivariant, namely the Hamiltonian flow $\Phi_t$ commutes with the action: $\Psi^{\hat{G}}_{(\sigma,\,g)}\circ \Phi_t= \Phi_t\circ \Psi^{\hat{G}}_{(\sigma,\,g)}\,\,\forall t\in\R,\forall (\sigma,\,g)\in \hat{G}$. 

Let $K<SO(3)$ be a finite subgroup of order $N$, with a fixed ordering  $K_o=(g_1=e,\,g_2,\,\dots,\,g_N)$. The subgroup $K$ acts linearly on $\mathbb{S}^2$, and we denote the set of points with non-trivial $K$-isotropy by:
\[\mathcal{F}[K]=\left\{r\in\mathbb{S}^2\,\colon \,K_r\ne \{e_K\}\right\}\]
where $K_r=\left\{g\in K\,\colon\, g\,r=r\right\}$ is the isotropy subgroup of $r$.
We introduce now the configurations of interest:
\begin{defn}[$K$-symmetric configuration] 
Denote by $M_{K_o}\subset M$ the subset of the phase space  defined by:
    \[M_{K_o}=\left\{\mathbf{r}\in M\,\colon\, r_i=g_i\,r_1\text{ for }i=1,\dots,N\right\}\,.\]
A configuration $\mathbf{r}\in M$ is said to be \textit{$K$-symmetric} if $\mathbf{r}\in M_{K_o}$ for some ordering $K_o$ of $K$.
\end{defn}
An important observation is that for a suitable discrete subgroup $\hat{K}<\hat{G}$ isomorphic to $K$, we have: \[M_{K_o}=\mathrm{Fix}(\hat{K})=\{\mathbf{r}\in M\,:\,\hat{g}\cdot\mathbf{r}=\mathbf{r}\,\,\,\forall \,\hat{g}\in\hat{K}  \}\,.\]In view of the $\hat{G}$-equivariance, this implies that $M_{K_o}$ is an invariant symplectic submanifold of the phase space and it makes sense to speak of $K$-symmetric \textit{solutions}. 

\subsection{The reduced system}

Keeping the same notation as above, introduce the map:
\begin{equation}\label{sympl diffeo}
\rho_{K_o}\colon \mathbb{S}^2\to M_{K_o}, \quad \rho_{K_o}(r)=(r,g_2\,r,\dots,g_N\,r)\,.
\end{equation}
This map defines a symplectic diffeomorphism between $(\mathbb{S}^2,\,N\omega_{\mathbb{S}^2})$ and $M_{K_o} $ endowed with (the restriction of) the form $\Omega$ in \eqref{form}. The solution of \eqref{system} with initial condition $\mathbf{r^0}\in M_{K_o}\setminus\mathscr{C}$ is $t\mapsto\rho_{K_o}(\tilde{r}(t))$, where $\tilde{r}(t)$ is the integral curve with initial condition $\tilde{r}(0)=r^0_1$ of the following \textit{reduced system} on $\mathbb{S}^2\setminus\mathcal{F}[K]$\footnote{Again, $\nabla h_{K}$ denotes the gradient (with respect to the standard Euclidean metric $(\,\R^{3},\,\cdot\,)$) of any extension of $h_{K}$ to an open set $U\subset\R^{3}$ containing $\mathbb{S}^2$.}: 

\begin{equation}
\dot{r}=-\frac{1}{N} r\times \nabla h_{K}(r)\,.\label{reddyn}
\end{equation} 
The function $h_{K}:\mathbb{S}^2\setminus\mathcal{F}[K]\to\R$ is the \textit{reduced Hamiltonian}, related to the Hamiltonian $H$ by: \begin{equation} h_{K}(r)\coloneqq H(\rho_{K_o}(r))=-\frac{N}{4}\sum_{j=2}^N \ln\big(\norm{r-g_j\, r}^2\big)\,.\label{redham}\end{equation} As the notation suggests, the reduced Hamiltonian $h_{K}$ is actually independent of the chosen ordering $K_o$, thanks to the $S_N$-invariance of the Hamiltonian $H$. The flow of the reduced system \eqref{reddyn} is complete on $\mathbb{S}^2\setminus\mathcal{F}[K]$. On the other hand, the reduced Hamiltonian $h_K$ is not defined on the set $\mathcal{F}[K]$: these are the points that correspond to collision configurations, as stated in the next proposition.

\begin{prop}\label{prop: n-uple collisions}  For any $r\in\mathcal{F}[K]$, the point $\rho_{K_o}(r)$ is a $K$-symmetric collision configuration whose collisions occur at the points of the orbit $K r$. Moreover, these are all $k$-tuple collisions, where $k=\abs{K_r}$.
\end{prop}
\begin{rem} We stress that the reduced system \eqref{reddyn} is not the result of a Poisson or a symplectic reduction. In fact, it is obtained by conjugating by $\rho_{K_
o}^{-1}$ the restriction to $M_{K_o}$ of the full system \eqref{system}. \end{rem} A relevant property of the reduced Hamiltonian is the following:

\begin{prop}[Invariance of the reduced Hamiltonian]\label{prop: invar red sys} 
 Let $N(K)$ be the normaliser of $K$ in $O(3)$ and suppose that a subgroup $L< O(3)$ satisfies $K\leq L\leq N(K)$. Then the reduced Hamiltonian \eqref{redham} is $L$-invariant.
\end{prop}

\begin{proof}
Since $K \trianglelefteq
L$, for any $l\in L$ the conjugation map $C_l\, \colon K\to K$ defined by $g\mapsto l^{-1}gl$ is an automorphism of $K$, therefore $\forall\, l\in L$: \begin{align*}h_{K}(r)&=-\frac{N}{4}\sum_{j=2}^N \ln\big(\norm{r-g_j\, r}^2\big)=-\frac{N}{4}\sum_{j=2}^N \ln\big(\norm{r-(l^{-1}g_j\,l)\, r}^2\big)=\\&=-\frac{N}{4}\sum_{j=2}^N \ln\big(\norm{l\,r-g_j\,(l\, r)}^2\big)=h_{K}(l\,r)\,.\end{align*}
\end{proof}

\paragraph{Regularisation of the Hamiltonian}

In order to avoid dealing separately with collisions, we work with a regularisation that includes the singularities in $\mathcal{F}[K]$. The \textit{regularised reduced Hamiltonian} is the following smooth \textit{polynomial} function $\tilde{h}_{K}\,\colon\,\mathbb{S}^2\to\R$: 
\begin{equation} \tilde{h}_{K}(r)\vcentcolon=\exp(-2 h_{K}(r))=\prod_{j=2}^N\norm{r-g_j\,r}^N\,,\label{regredham}
\end{equation}
and the \textit{regularised reduced system} is given by the following smooth Hamiltonian vector field on $(\mathbb{S}^2,\,N\,\omega_{\,\mathbb{S}^2})$:
\begin{equation}
    \dot{r}=-\frac{1}{N}r\times\nabla\tilde{h}_{K}(r)\label{regredsys}
\end{equation} 
This modification of the Hamiltonian defines a complete flow on $\mathbb{S}^2$ by reparametrising the orbits of the reduced system \eqref{reddyn} and artificially introducing the points in $\mathcal{F}[K]$ as stable equilibria. Motivated by this observation, the points in $\mathcal{F}[K]$ will be called \textit{collision equilibria} of the regularised reduced system \eqref{regredsys}, while the other equilibria of \eqref{regredsys} will be called \textit{non-collision equilibria}. Notice that collision equilibria are the only points in the $0$ level set of \eqref{regredham}.
Furthermore, the level sets in $\mathbb{S}^2\setminus\mathcal{F}[K]$ of the Hamiltonians $h_K$ and $\tilde{h}_K$ coincide (even if they correspond to different level values). 
\begin{rem} Observe that both the reduced system \eqref{reddyn} and the regularised reduced system \eqref{regredsys} are integrable 1-degree-of-freedom Hamiltonian systems.
\end{rem}

\section{Finding critical points of the regularised Hamiltonian}
\begin{namedthm*}{Notation}
For the rest of the article, the full octahedral group $\mathbb{O}_h$ is identified with the subgroup of $O(3)$ generated by the real reflections that send the octahedron with the standard orientation\footnote{Vertices: $(\pm 1,0,0)$, $(0,\pm 1,0)$, $(0,0,\pm 1)$.} to itself. Similarly,
the full icosahedral group $\mathbb{I}_h$ is identified with the subgroup of $O(3)$ generated by the real reflections that send the icosahedron with the standard orientation\footnote{Vertices: even permutations of $\frac{1}{\sqrt{1+\phi^2}}(\pm 1,\pm\phi,0)$.} to itself. Accordingly, we regard the octahedral group as $\mathbb{O}=\mathbb{O}_h\cap SO(3)$ and the icosahedral group as $\mathbb{I}=\mathbb{I}_h\cap SO(3)$.
\end{namedthm*}
In view of the previous section, our goal becomes classifying the critical points of the regularised reduced Hamiltonian:
\begin{equation}\label{regredham2}\tilde{h}_{K}\,\colon\,\mathbb{S}^2\to\R,\quad \tilde{h}_{K}(r)=\prod_{\substack{g\in K\\g\ne e}}\norm{r-g\,r}^{N}\,, \quad \text{where }K=\mathbb{O}\text{ or }K=\mathbb{I}\,.\end{equation} Let $K_h$ be the full octahedral group $\mathbb{O}_h<O(3)$ or the full icosahedral group $\mathbb{I}_h<O(3)$, respectively. Observe that Proposition \ref{prop: invar red sys} implies that $\tilde{h}_{K}$ is $K_h$-invariant. Our strategy relies on the $K_h$-invariance of $\tilde{h}_{K}$ to find a chart of $\mathbb{S}^2$ that simplifies the analysis, proceeding as follows.

\subsection{Fundamental region}
\begin{figure}
    \centering
    \subfloat[Tessellation of $\mathbb{S}^2$ induced by $\mathbb{O}_h$.]{\includegraphics[width=0.45\linewidth]{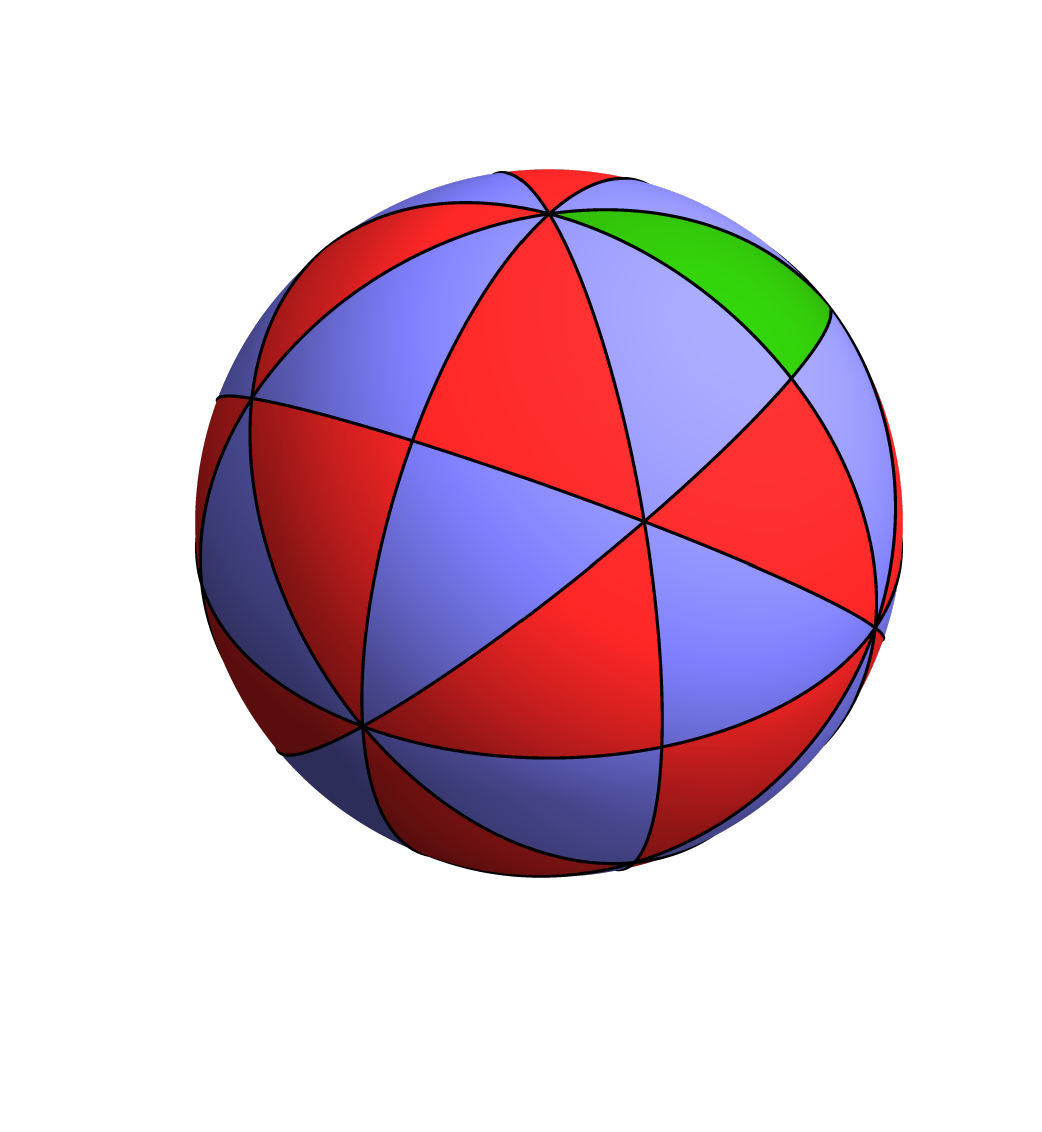}}
        \quad
   \subfloat[Tessellation of $\mathbb{S}^2$ induced by $\mathbb{I}_h$.]{\includegraphics[width=0.45\linewidth]{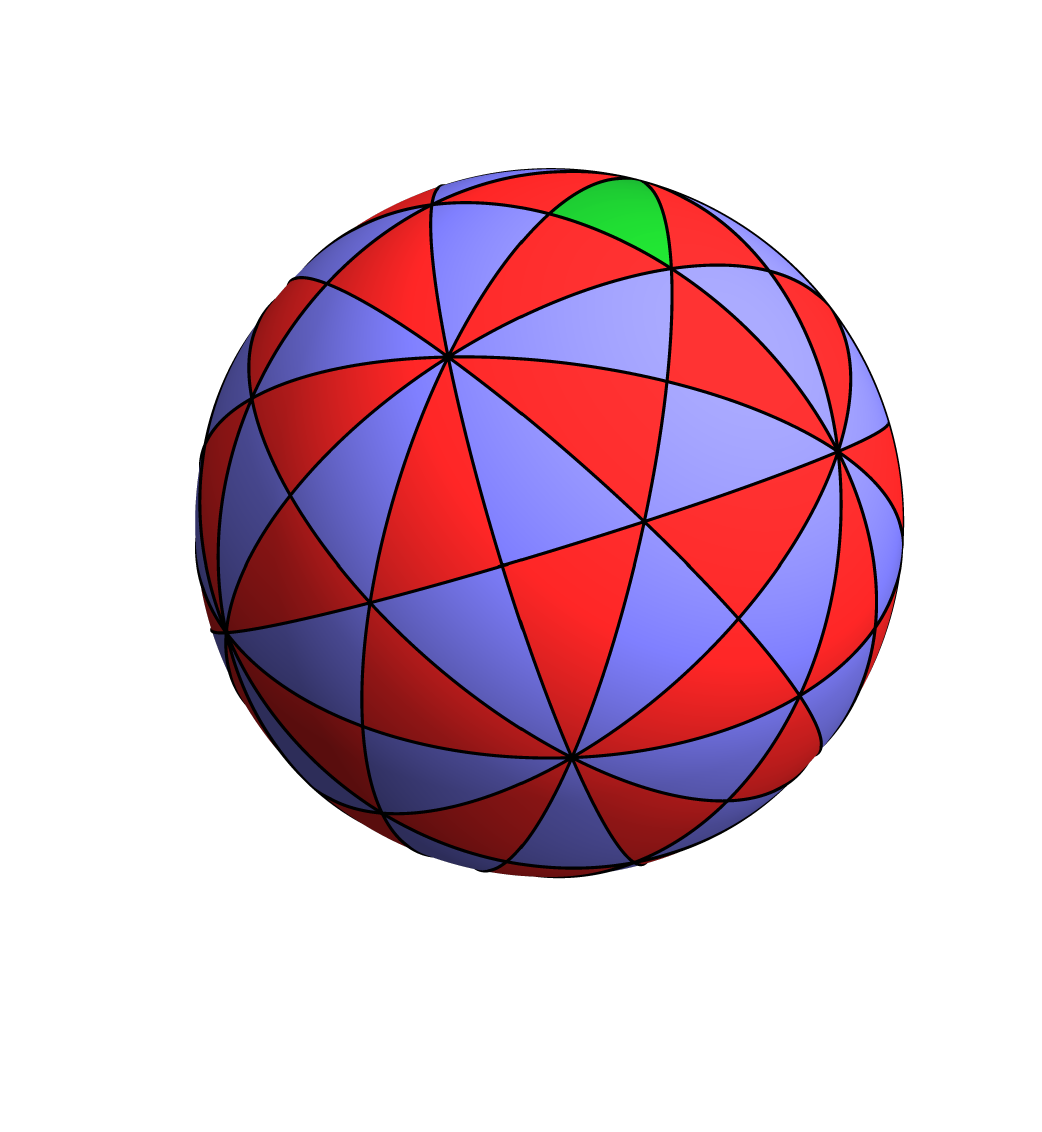}}
    
    \caption{A picture of the tessellation of $\mathbb{S}^2$ produced by the action of $\mathbb{O}_h$ and $\mathbb{I}_h$. The green spherical triangles (together with their edges and vertices) are the fundamental domains $\mathcal{R}\subset\mathbb{S}^2$ we will restrict to.}\label{fig: Tess}
\end{figure}

First, we identify a convenient fundamental region $\mathcal{R}\subset\mathbb{S}^2$ for the linear action of $K_h$ on $\mathbb{S}^2$. We can then restrict the analysis of $\tilde{h}_{K}$ to $\mathcal{R}$ due to $K_h$-invariance. 
The theory of real reflection groups prescribes a method to obtain the tessellation of $\mathbb{S}^2$ into fundamental triangles: it is enough to consider the symmetry planes of the octahedron/icosahedron. The sphere is tiled by such planes into $\abs{K_h}=2N$ fundamental triangles, and we can choose the most suitable (see Figure \ref{fig: Tess}). Each fundamental triangle is the intersection of a Weyl chamber with the sphere.
    
Observe that the only points with non-trivial $K_h$-isotropy are those situated on an edge of one of the fundamental triangles. Moreover, we have that $\mathcal{F}[K]$ coincides with the vertices of the fundamental triangles.

\subsection{Parametrising the fundamental region} 

Once the analysis is restricted to a fundamental region $\mathcal{R}$, we rely  on the following classic theorem (see \cite{Hum}) to find a parametrisation of the interior $\mathrm{int}(\mathcal{R})$ that significantly lowers the degree of the regularised Hamiltonian \eqref{regredham2}. 

\begin{thm}[Chevalley--Shephard--Todd Theorem]\label{thm: CST}
    Let $G<GL_n(\R)$ be a finite matrix group. The algebra of $G$-invariant polynomials in $n$ variables $\R[x_1,\dots,x_n]^G$ is a polynomial algebra generated by $n$ homogeneous algebraically independent invariants $\theta_1,\dots,\theta_n$ if and only if $G$ is a finite reflection group (i.e. it is generated by reflections).
\end{thm}
 Applied to the current case, the theorem ensures the existence of three homogeneous algebraically independent invariants $\theta_1$, $\theta_2$ and $\theta_3$ (clearly, those invariants differ in the two cases $K_h=\mathbb{O}_h$ and $K_h=\mathbb{I}_h$) such that \[\R[x_1,x_2,x_3]^{K_{h}}=\R[\theta_1,\theta_2,\theta_3]\, .\] In general, the invariants are not uniquely determined by the group $G$, yet their degrees are. If their degrees are known, a set of invariants can be built by applying the \textit{Reynolds operator} (see \cite{Sturm}) to a set of polynomials of the corresponding degrees and checking that the resulting invariants are actually algebraically independent. We choose the following invariants in the case $K_h=\mathbb{O}_h$:
\begin{equation}\label{oct inv}
\begin{aligned}
    &\theta_1(x,y,z)\coloneqq x^2+y^2+z^2\,,\\
    &\theta_2(x,y,z)\coloneqq x^2y^2+x^2z^2+y^2z^2\,,\\
    &\theta_3(x,y,z)\coloneqq x^2y^2z^2\,.
\end{aligned}
\end{equation}
In the case $K_h=\mathbb{I}_h$, we take again $\theta_1(x,y,z)\coloneqq x^2+y^2+z^2$, while the complete expressions of the other two invariants are given in Appendix B (see \eqref{icoinv6}, \eqref{icoinv10}).

Denote by $\mathscr{I}\,\colon\,\R^3\to\R^3$ the map defined by such invariants as follows: \[\mathscr{I}(x,y,z)\coloneqq(\theta_1(x,y,z),\,\theta_2(x,y,z),\,\theta_3(x,y,z))\,.\]

With our choice, the second and third fundamental invariants $\theta_2$ and $\theta_3$ provide local coordinates for the interior of the fundamental region $\mathcal{R}$. This is a consequence of the following Proposition, whose proof can be found in \cite{Poma}:

\begin{prop}\label{prop: invariants as parameters} The restriction $\restr{\mathscr{I}}{\mathrm{int}(\mathcal{R})}$ to $\mathrm{int}(\mathcal{R})\subset\mathbb{S}^2$ is a diffeomorphism onto its image.
\end{prop}

Dropping the first entry which is identically equal to $1$, the image $\mathscr{I}(\mathrm{int}(\mathcal{R}))$ can be identified with the open set \[U\coloneqq\{(\theta_2(x,y,z),\,\theta_3(x,y,z)):(x,y,z)\in\mathrm{int}(\mathcal{R})\}\subset \R^2\,.\] In view of Proposition \ref{prop: invariants as parameters}, we introduce the following chart of $\mathrm{int}(\mathcal{R})$:
\[\psi^{-1}\,\colon \,\mathrm{int}(\mathcal{R})\to U\quad (x,y,z)\mapsto(\theta_2(x,y,z),\,\theta_3(x,y,z))\,,\]
whose inverse defines the parametrisation:
\begin{equation}\label{param}
    \psi\,\colon \,U\to\mathrm{int}(\mathcal{R})\quad(s,t)\mapsto \psi(s,t)\,.
\end{equation} Hence, the task of finding the critical points of $\tilde{h}_K$ in $\mathrm{int}(\mathcal{R})$ reduces to the easier task of finding the critical points in $U$ of the polynomial in $(s,t)$ given by: \[\tilde{h}_K^I(s,t)\coloneqq \tilde{h}_K(\psi(s,t))\,\,\,\text{for}\,\,(s,t)\in U\,.\] Observe that the degree of $\tilde{h}_K^I$ is considerably lower than the degree of $\tilde{h}_K$: an upper bound is $\deg \tilde{h}_K^I\le \big\lfloor\frac{\deg{\tilde{h}_K}}{\deg{\theta_2}}\big\rfloor$. 

\subsection{The boundary of the fundamental region}
\begin{figure}
    \centering
    \subfloat[The chosen fundamental region $\mathcal{R}$ in the case $K_h=\mathbb{O}_h$.\label{fig: Full Oct Fund Reg}]{\raisebox{0.05\height}{\begin{overpic}[width=0.4\textwidth]{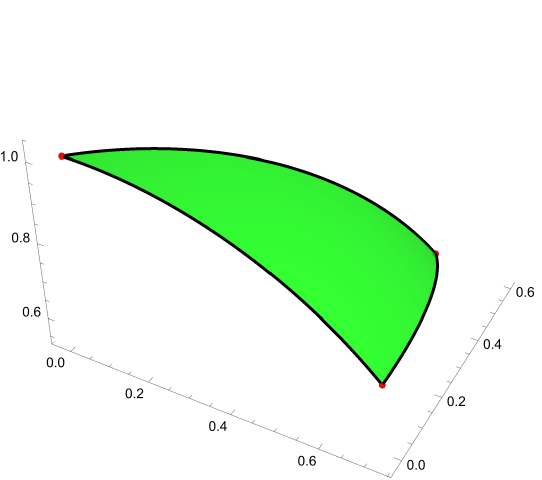}
    \put(30,8){$x$}
    \put(90,10){$y$}
    \put(0,40){$z$}
    \put(8,63){\textcolor{red}{$V_1$}}
    \put(70,10){\textcolor{red}{$V_2$}}
    \put(82,40){\textcolor{red}{$V_3$}}
    \put(33,40){$\mathcal{E}_{12}$} 
    \put(58,58){$\mathcal{E}_{13}$}
    \put(80,30){$\mathcal{E}_{23}$}
    \end{overpic}}}
    \quad
    \subfloat[The chosen fundamental region $\mathcal{R}$ in the case $K_h=\mathbb{I}_h$.\label{fig: Full Ico Fund Reg}]{\raisebox{0.1\height}{\begin{overpic}[width=0.45\textwidth]{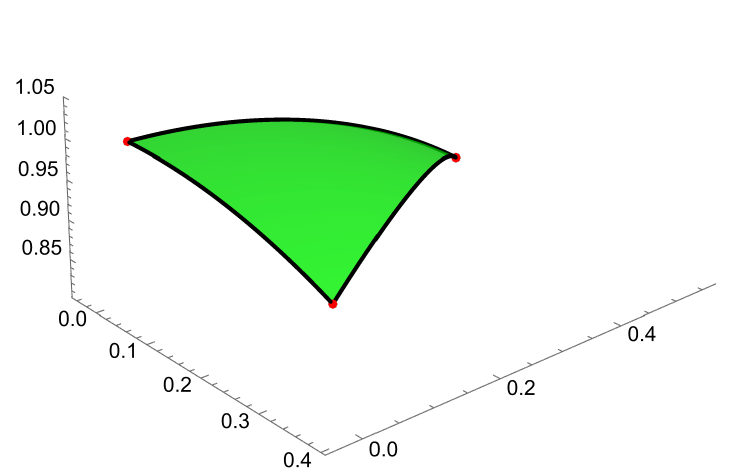}
    \put(20,8){$x$}
    \put(75,10){$y$}
    \put(-2,40){$z$}
    \put(11,47){\textcolor{red}{$V_1$}}
    \put(43,16){\textcolor{red}{$V_2$}}
    \put(63,40){\textcolor{red}{$V_3$}}
    \put(22,32){$\mathcal{E}_{12}$} 
    \put(38,50){$\mathcal{E}_{13}$}
    \put(55,30){$\mathcal{E}_{23}$}
    \end{overpic}}}
    
    \caption{A picture of the chosen fundamental regions $\mathcal{R}\subset\mathbb{S}^2$ with their vertices and edges.}\label{fig: Fund Reg}
\end{figure}

\begin{namedthm*}{Notation} For the rest of the article, for $\alpha\in \R^3\setminus\{0\}$, we denote by $r_{\alpha}\in O(3)$ the (real) reflection with respect to the plane $H_\alpha$ orthogonal to $\langle\alpha\rangle$.\end{namedthm*} The boundary $\partial\mathcal{R}$ of the fundamental region consists of three circular arcs (\textit{edges}) that meet only at the vertices $V_1$, $V_2$ and $V_3$. Define $\mathcal{E}_{ij}$ as the edge with vertices $V_i$ and $V_j$, so that $\mathcal{E}_{12}\cap\mathcal{E}_{13}=V_1$, $\mathcal{E}_{12}\cap\mathcal{E}_{23}=V_2$, $\mathcal{E}_{13}\cap\mathcal{E}_{23}=V_3$ and  $\partial\mathcal{R}=\mathcal{E}_{12}\cup\mathcal{E}_{13}\cup\mathcal{E}_{23}$ (see Figure \ref{fig: Fund Reg}). Each edge $\mathcal{E}_{ij}$ lies on a great circle $\mathbb{S}^2\cap H_{\alpha_{ij}}$, where $\alpha_{ij}$ is a suitably chosen vector with the property that $r_{\alpha_{ij}}\in K_h$.

\paragraph{Critical points on the edges.}\label{subsec: edges fund reg}
The map $\mathscr{I}$ defined above is not a local diffeomorphism at any point $r\in\partial \mathcal{R}$, as it fails to be injective in any neighbourhood containing $r$. Hence, it cannot be used to define a local chart to study critical points on the edges. Instead of working with a chart at $r$ which is unrelated to the $K_h$-symmetry, we invoke the following principle (see \cite{Pal}):

\begin{thm}[Palais' principle of symmetric criticality] Let $G$ be a group of isometries of a Riemannian manifold $X$ and let $f\in C^1(X)$ be a function invariant under $G$. Then the set $\Sigma$ of fixed points of $X$
under the action of $G$ is a (totally geodesic) smooth submanifold of $X$, and if $p\in\Sigma$ is a
critical point of $\restr{f}{\Sigma}$ then $p$ is in fact a critical point of $f$.
\end{thm}

The principle tells us that if we are interested in finding the critical points of any $K_h$-invariant function $f\in C^{1}(\mathbb{S}^2)$ on an edge of $\mathcal{R}$, we can restrict it to the $1$-dimensional submanifold of $\mathbb{S}^2$ consisting of the great circle $\Sigma$ containing the edge. 
In fact, this corresponds to taking $X=\mathbb{S}^2$ and $G=\langle \,r_\alpha\,\rangle\leq K_h$ in the statement of Palais' principle, where $r_\alpha$ is the reflection with respect to the hyperplane $H_\alpha$ containing the edge. We have then $\Sigma=H_\alpha\cap\mathbb{S}^2$. This reduces the problem to finding the critical points of a smooth function of one variable defined on the edges $\mathcal{E}_{12},\mathcal{E}_{13}, \mathcal{E}_{23}\subset\partial\mathcal{R}$.

\paragraph{The vertices are always critical points.}
Let $V\in\partial\mathcal{R}$ be any of the three vertices and let $\Sigma_i=H_{\alpha_i}\cap\mathbb{S}^2$, $i=1,2$, be the two great circles containing the edges of $\mathcal{R}$ adjacent to $V$. Since $V\in H_{\alpha_i}$, we have $\alpha_i\in T_V\mathbb{S}^2$ and $d_Vr_{\alpha_i}(\alpha_i)=r_{\alpha_i}(\alpha_i)=-\alpha_i$; hence, for any $K_h$-invariant $f\in C^1(\mathbb{S}^2)$, differentiating $f\circ r_{\alpha_i}=f$ at $V$ gives $d_Vf(\alpha_i)=d_Vf(-\alpha_i)$, i.e. $d_Vf(\alpha_i)=0$. As $H_{\alpha_1}\neq H_{\alpha_2}$, the vectors $\alpha_1,\alpha_2$ span $T_V\mathbb{S}^2$, so $d_Vf=0$ and $V$ is a critical point of $f$.

\section{Octahedrally symmetric solutions of the $24$-vortex problem on the sphere}

We apply the method of the previous section to the case $K=\mathbb{O}$ and $N=\abs{K}=24$. For simplicity, rather than working with the regularised reduced Hamiltonian $\tilde{h}_{\mathbb{O}}$, we can equivalently consider the polynomial function \[p_{\mathbb{O}}(x,y,z)\vcentcolon=\frac{1}{2^{26}}\sqrt[12]{\tilde{h}_{\mathbb{O}}(x,y,z)}\,,\]
which has the same critical points and the same level sets as $\tilde{h}_{\mathbb{O}}$ (corresponding to different values, though). The complete expressions of $\tilde{h}_{\mathbb{O}}$ and $p_{\mathbb{O}}$ are given in Appendix B (see \eqref{octregham}, \eqref{octsimplham}). The full octahedral group $\mathbb{O}_h$ is generated by the following reflections:
\begin{align*}
    &r_{e_y}=\begin{pmatrix}
1 & 0 &0\\
0 & -1&0\\
0&0&1
\end{pmatrix}\,,\\
&r_{e_x-e_y}=\begin{pmatrix}
0 & 1 &0\\
1 & 0&0\\
0&0&1
\end{pmatrix}\,,\\
&r_{-e_x+e_z}=\begin{pmatrix}
0 & 0 &1\\
0 & 1&0\\
1&0&0
\end{pmatrix}\,.
\end{align*}
The reflection planes of such generators contain the three edges of the fundamental region $\mathcal{R}$ in green in Figure \ref{fig: Full Oct Fund Reg}. 
Explicitly, $\mathcal{R}$ is described by:
\[\mathcal{R}=\{(x,y,z)\in\mathbb{S}^2:0\le y\le x\le z\}\,.\] We denote by $\mathcal{E}_{12}$ the edge $\mathcal{R}\cap H_{e_y}$, by $\mathcal{E}_{13}$ the edge $\mathcal{R}\cap H_{e_x-e_y}$ and by  $\mathcal{E}_{23}$ the edge $\mathcal{R}\cap H_{-e_x+e_z}$. The three vertices of $\mathcal{R}$ are:
\[\begin{aligned} &V_1\coloneqq\mathcal{E}_{12}\cap\mathcal{E}_{13}=(0,0,1)\,,\\ &V_2\coloneqq\mathcal{E}_{12}\cap\mathcal{E}_{23}=(1/\sqrt{2},0,1/\sqrt{2})\,,\\ &V_3\coloneqq\mathcal{E}_{13}\cap\mathcal{E}_{23}=(1/\sqrt{3},1/\sqrt{3},1/\sqrt{3})\,.\end{aligned}\] Finally, we describe the isotropy group of points on $\partial\mathcal{R}$: these are the only points with non-trivial isotropy group. The isotropy group of any point $r\in\mathcal{E}_{12}\setminus\{V_1,V_2\}$ is $(\mathbb{O}_h)_r=\langle \,r_{e_y}\,\rangle$; similarly, for any point $r\in\mathcal{E}_{13}\setminus\{V_1,V_3\}$ we have $(\mathbb{O}_h)_r=\langle \,r_{e_x-e_y}\,\rangle$ and for any point $r\in\mathcal{E}_{23}\setminus\{V_2,V_3\}$ we have $(\mathbb{O}_h)_r=\langle\, r_{-e_x+e_z}\,\rangle$. As for the vertices, $(\mathbb{O}_h)_{V_1}=\langle\, r_{e_y},\,r_{e_x-e_y}\,\rangle$, $(\mathbb{O}_h)_{V_2}=\langle\, r_{e_y},\,r_{-e_x+e_z}\,\rangle$ and $(\mathbb{O}_h)_{V_3}=\langle\, r_{e_x-e_y},\,r_{-e_x+e_z}\,\rangle$.

\subsection{The polynomial $p_{\mathbb{O}}$ in invariant coordinates}\label{sec: oct eqs}

Before dealing with the classification of critical points, we rename the invariants in \eqref{oct inv} as follows: 
\begin{align*}&s\coloneqq \theta_2(x,y,z)= x^2y^2+x^2z^2+y^2z^2\,,\\ &t\coloneqq \theta_3(x,y,z)=x^2y^2z^2\,.\end{align*}
Using the parametrisation $\psi$ in \eqref{param}, we obtain the polynomial function $p^I_{\mathbb{O}}\,\colon U\to\R$ in the variables $(s,t)$.
An explicit expression of $p^I_{\mathbb{O}}$ is the following (for details about its derivation, see Remark 3.2.4 of \cite{Poma}):

\begin{equation}\label{reg red oct ham in the inv}\begin{aligned} &p^I_{\mathbb{O}}(s,t)=F_1(s,t)^2 F_2(s,t)^3 F_3(s,t),\,\text{ where:}\\
   &F_1(s,t)\coloneqq 1 - 2 s + s^2 - 12 t\,,\\
   &F_2(s,t)\coloneqq t - s\,,\\
   &F_3(s,t)\coloneqq -4 + 15 s^2 + 4 s^3 - 144 t + 
   27 t^2 + 6 s (2 + 15 t)\,.
   \end{aligned}\end{equation}

\begin{figure}
   \centering
   \begin{overpic}[width=0.98\textwidth]{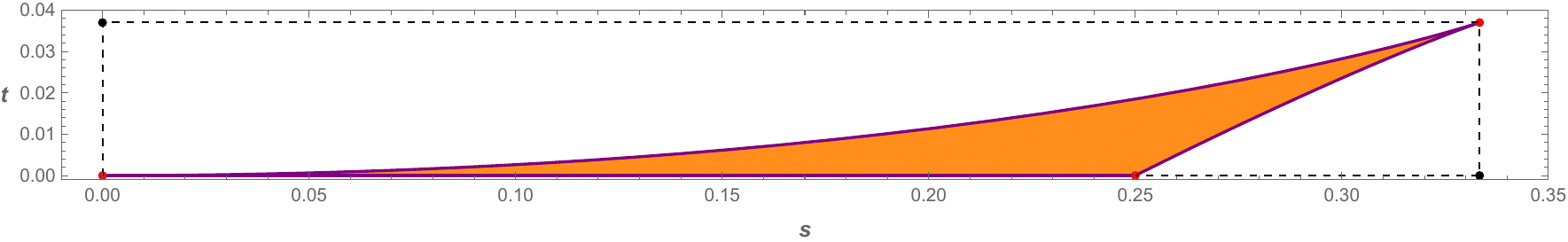}
    \put(7,5){\textcolor{red}{$V_1$}}
    \put(69,5){\textcolor{red}{$V_2$}}
    \put(95,12){\textcolor{red}{$V_3$}}
    \put(88,1){\textcolor{black}{$(1/3,\,0)$}}
    \put(7,11){\textcolor{black}{$(0,\,1/27)$}}
    \end{overpic}
    \caption{In orange, the domain $U=\psi^{-1}(\mathrm{int}(\mathcal{R}))$ in the invariants $(s, t)$-plane that parametrises the interior $\mathrm{int}(\mathcal{R})$ of the fundamental region.}
    \label{fig: FullOctInvRegNoGH}
\end{figure}

\begin{rem}\label{rem: oct inv fund reg bounds} Let $(x,y,z)$ be any point in $\mathrm{int}(\mathcal{R})$ and $(s,t)=\psi^{-1}(x,y,z)$ be the corresponding point in $U$. Then, the conditions $x^2+y^2+z^2=1$ and $0<y<x<z$ coming from the definition of $\mathrm{int}(\mathcal{R})$ impose the following bounds: \[0<\sqrt[3]{t}=\sqrt[3]{x^2\,y^2\,z^2}< \frac{x^2+y^2+z^2}{3}=\frac{1}{3}\,,\]
\[0<s=x^2y^2+x^2z^2+y^2z^2< x^4+y^4+z^4=(x^2+y^2+z^2)^2-2(x^2y^2+x^2z^2+y^2z^2)=1-2s\,,\]
showing that $U\subseteq (0,\,1/3)\times (0,\,1/27)$ (see Figure \ref{fig: FullOctInvRegNoGH}).
\end{rem}

\subsection{Classification of octahedrally symmetric equilibrium and collision configurations}
\begin{figure}
    \centering
    \subfloat[Pseudo snub cube equilibrium configuration.\label{fig: oct non collision eqs (a)}]{\includegraphics[width=0.4\linewidth]{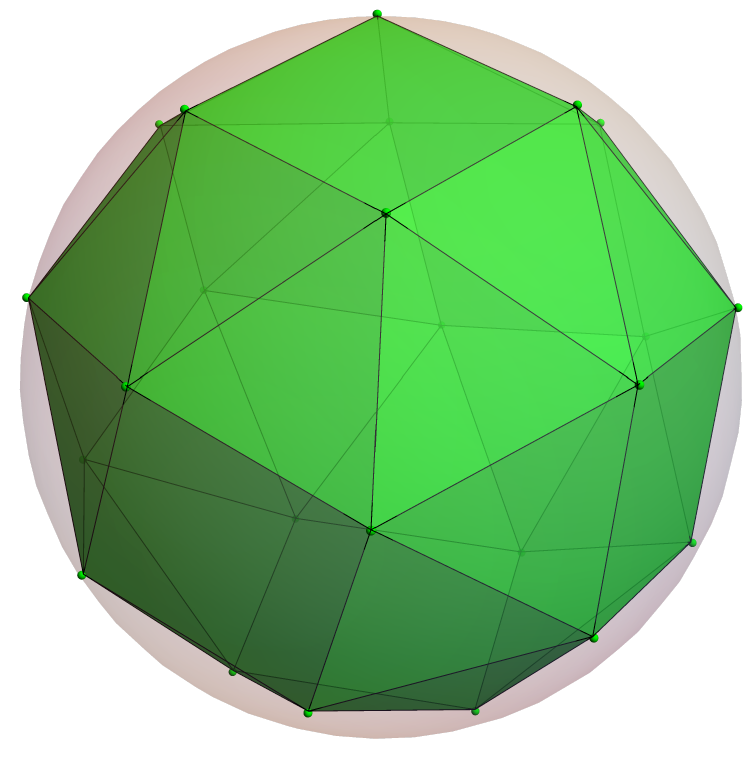}}
        \quad
    \subfloat[Truncated octahedron equilibrium configuration.\label{fig: oct non collision eqs (b)}]{\includegraphics[width=0.4\linewidth]{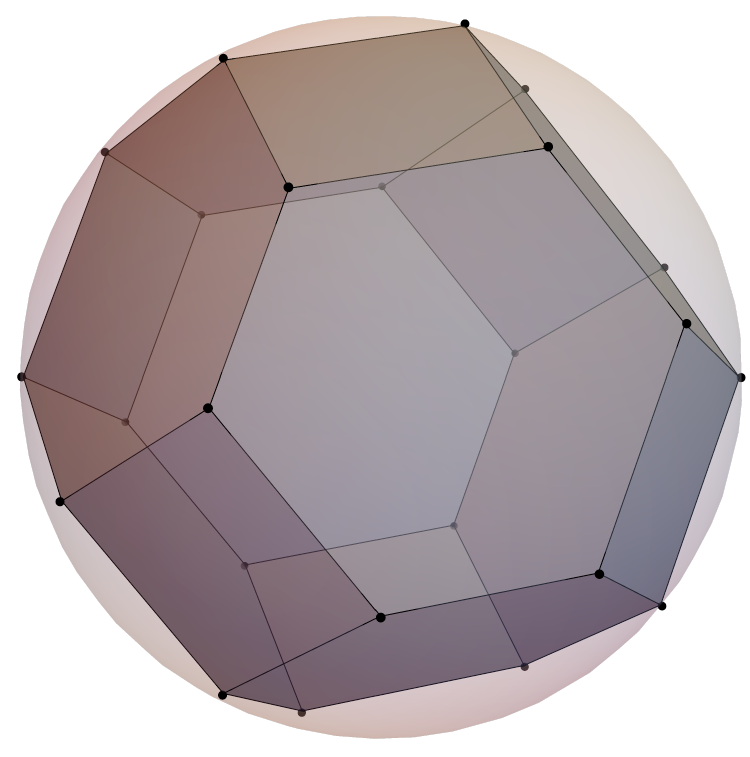}}
    \quad
    \subfloat[Pseudo rhombicuboctahedron equilibrium configuration.\label{fig: oct non collision eqs (c)}]{\includegraphics[width=0.4\linewidth]{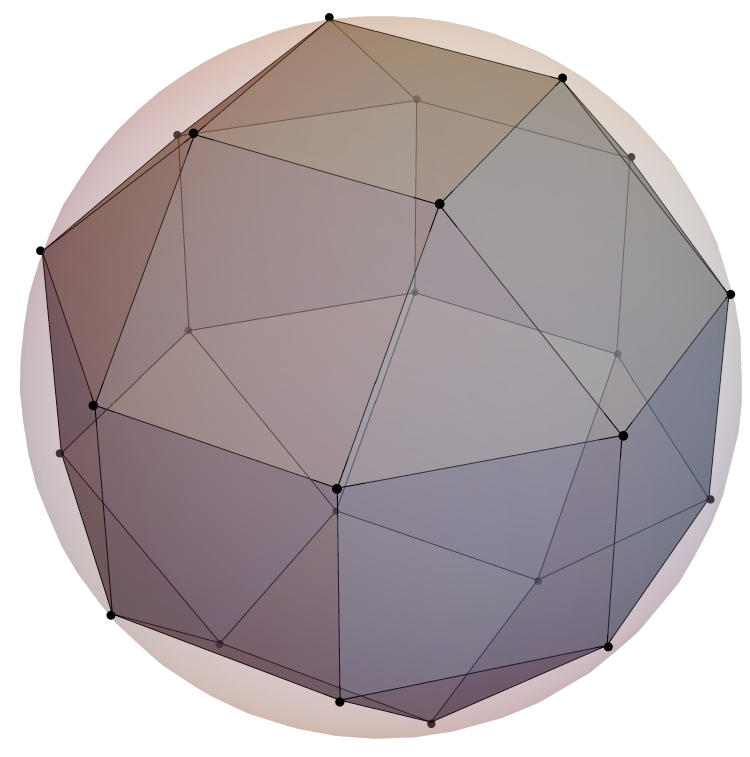}}
    \quad
    \subfloat[Truncated cube equilibrium configuration.\label{fig: oct non collision eqs (d)}]{\includegraphics[width=0.4\linewidth]{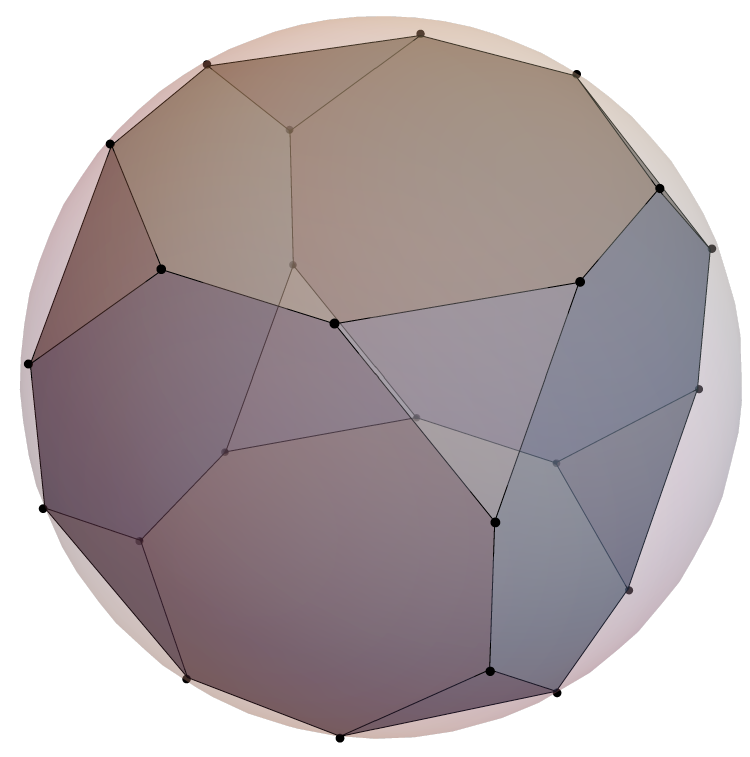}}
    \caption{Equilibrium configurations of the full system \eqref{system} described in item (i) of Theorem \ref{thm: Oct eqs}. For more details on the dimensions, see Appendix A.}\label{fig: oct non collision eqs}
\end{figure}

\begin{figure}
    \centering
    \subfloat[Octahedral collision configuration: simultaneous quaternary collision at each vertex.\label{fig: oct collision eqs (a)}]{\includegraphics[width=0.4\linewidth]{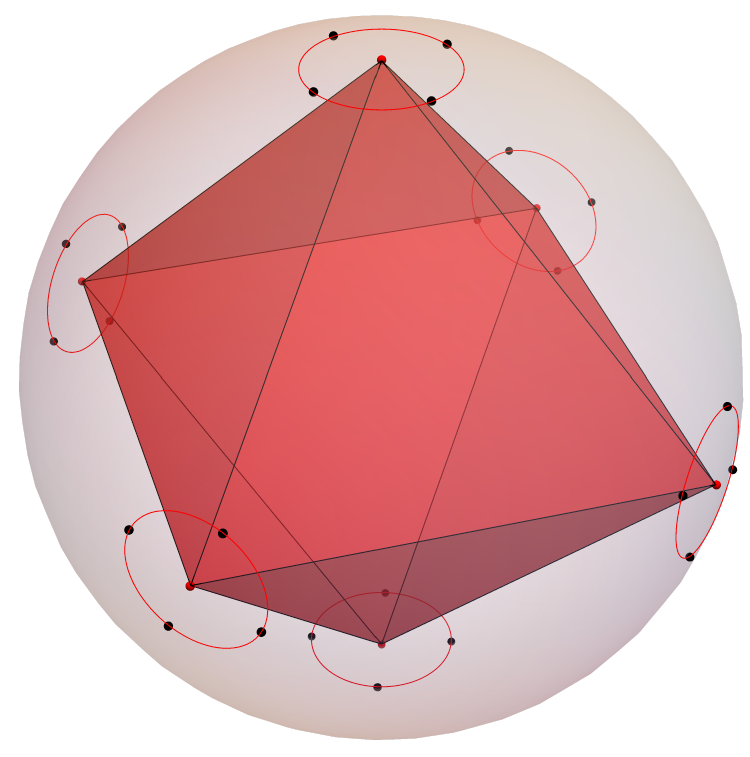}}
    \quad
    \subfloat[Cuboctahedral collision configuration: simultaneous binary collision at each vertex.\label{fig: oct collision eqs (b)}]{\includegraphics[width=0.4\linewidth]{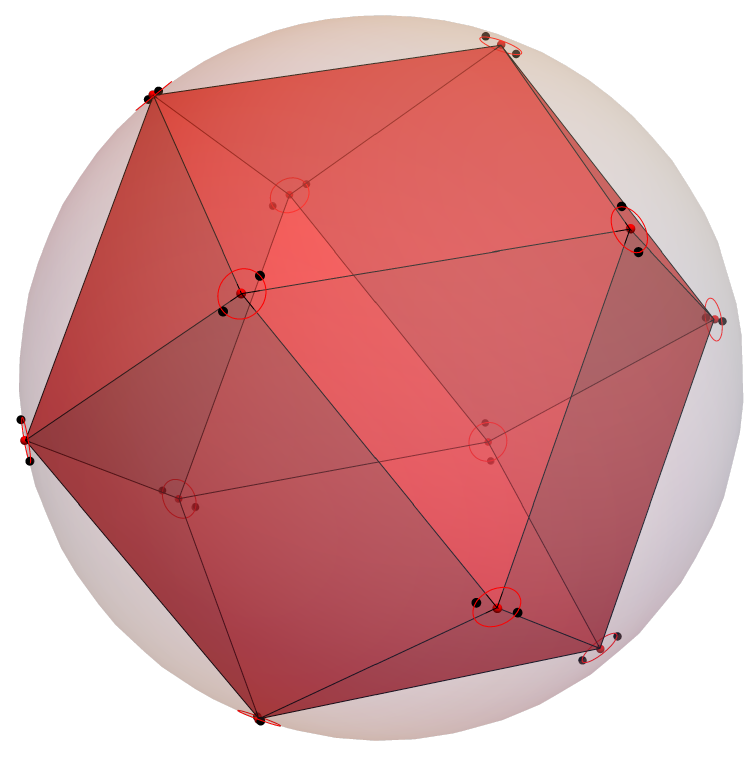}}
    \quad
    \subfloat[Cubic collision configuration: simultaneous ternary collision at each vertex.\label{fig: oct collision eqs (c)}]{\includegraphics[width=0.4\linewidth]{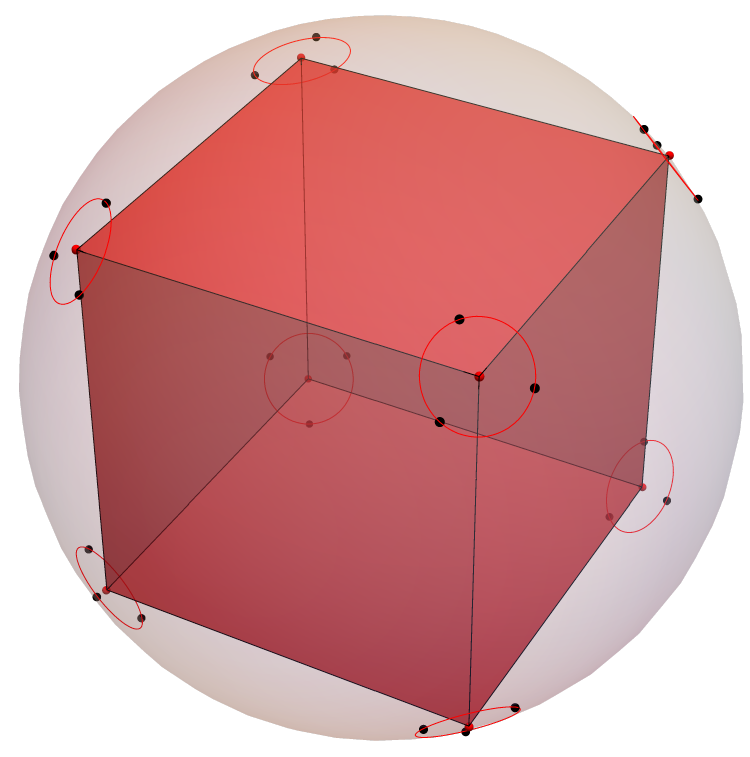}}
    \caption{Collision configurations of the full system \eqref{system} described in item (ii) of Theorem \ref{thm: Oct eqs}.} \label{fig: oct collision eqs}
\end{figure}

Our main classification result in this case is the following:
\begin{thm}[Octahedrally symmetric equilibrium and collision configurations] \label{thm: Oct eqs}  Let $K=\mathbb{O}$, and $N=24$. 
\begin{itemize}
\item[(i)]The following are the only non-collision equilibria of the regularised reduced system \eqref{regredsys}:
\begin{itemize}

    \item[(a)] The \textbf{pseudo}\footnote{We use the word “pseudo” to point out that it has the same vertex configuration but irregular faces.} \textbf{snub cube equilibria} occurring at the two distinct $\mathbb{O}$-orbits contained in the $\mathbb{O}_h$ orbit of the point: \[(\alpha_1,\alpha_2,\alpha_3)\in \mathrm{int}(\mathcal{R})\]
    whose coordinates $0<\alpha_2<\alpha_1<\alpha_3<1$ satisfy: \[y^3-y^2+\bar{s}y-\bar{t}=(y-\alpha_1^2)(y-\alpha_2^2)(y-\alpha_3^2)\] where $\bar{s}$ is the only root in the interval $(0,\,1/3)$ of a polynomial $q_1(s)$ of degree $8$ and $\bar{t}$ is the only root in the interval $(0,\,1/27 )$ of a polynomial $q_2(t)$ of degree $8$ (see \eqref{octinteqpoly1}, \eqref{octinteqpoly2} for their complete expression). These are stable equilibria of the regularised reduced system \eqref{regredsys} which correspond to equilibrium configurations of the full system \eqref{system} where the vortices are situated at the vertices of an $\mathbb{S}^2$-inscribed left-chiral pseudo snub cube or its right-chiral version (see Figure \ref{fig: oct non collision eqs (a)} for the right-chiral version). The point $(\alpha_1,\alpha_2,\alpha_3)$ defined above is approximately given by $(0.439584, 0.250806, 0.862474)$.
    
    \item[(b)] The \textbf{truncated octahedron equilibria} occurring at the $\mathbb{O}$-orbit of the point:
    \[(\sqrt{\beta},0,\sqrt{1-\beta})\in \mathcal{E}_{12}\]
    where $\beta$ is the only root of $p_1(u)\coloneqq -6+41 u-72 u^2+22 u^3+89 u^4-120 u^5+40 u^6$ in the interval $(0,\,1/2)$. These are unstable equilibria of the regularised reduced system \eqref{regredsys} which correspond to equilibrium configurations of the full system \eqref{system} where the vortices are situated at the vertices of an $\mathbb{S}^2$-inscribed irregular truncated octahedron (see Figure \ref{fig: oct non collision eqs (b)}). The point $(\sqrt{\beta},0,\sqrt{1-\beta})$ defined above is approximately given by $(0.474726, 0, 0.880134)$.
    
    \item[(c)] The \textbf{pseudo rhombicuboctahedron equilibria} occurring at the $\mathbb{O}$-orbit of the point:
    \[(\sqrt{\gamma},\sqrt{\gamma},\sqrt{1-2\gamma})\in \mathcal{E}_{13}\]
    where $\gamma$ is the only root of $p_2(u)\coloneqq -3+41 u-190 u^2+502 u^3-323 u^4-1023 u^5+1188 u^6$ in the interval $(0,1/3)$. These are unstable equilibria of the regularised reduced system \eqref{regredsys} which correspond to equilibrium configurations of the full system \eqref{system} where the vortices are situated at the vertices of an $\mathbb{S}^2$-inscribed pseudo rhombicuboctahedron (see Figure \ref{fig: oct non collision eqs (c)}). The point $(\sqrt{\gamma},\sqrt{\gamma},\sqrt{1-2\gamma})$ defined above is approximately given by $(0.350294, 0.350294, 0.868671)$.

    \item[(d)] The \textbf{truncated cube equilibria} occurring at the $\mathbb{O}$-orbit of the point:
    \[(\sqrt{\delta},\sqrt{1-2\delta},\sqrt{\delta})\in \mathcal{E}_{23}\]
    where $\delta$ is the only root of $p_2(u)$ in the interval $(1/3,\,1/2)$. These are unstable equilibria of the regularised reduced system \eqref{regredsys} which correspond to equilibrium configurations of the full system \eqref{system} where the vortices are situated at the vertices of an $\mathbb{S}^2$-inscribed irregular truncated cube (see Figure \ref{fig: oct non collision eqs (d)}). The point $(\sqrt{\delta},\sqrt{1-2\delta},\sqrt{\delta})$ defined above is approximately given by $(0.685580, 0.244870, 0.685580)$.
\end{itemize}

\item[(ii)]The only collision equilibria of the regularised reduced system \eqref{regredsys} are:

\begin{itemize}

\item[(a)] The \textbf{octahedral collision equilibria} occurring at the vertices of an octahedron $(\pm 1, 0,0)$, $(0,\pm 1,0)$, $(0,0,\pm 1)$, which form the $\mathbb{O}$-orbit of the vertex $V_1\in\mathcal{R}$. These are stable equilibria of the regularised reduced system \eqref{regredsys} which correspond to collision configurations of the full system \eqref{system} where six quadruplets of vortices simultaneously collide at the vertices of an $\mathbb{S}^2$-inscribed regular octahedron (see Figure \ref{fig: oct collision eqs (a)}).

\item[(b)] The \textbf{cuboctahedral collision equilibria} occurring at the vertices of a cuboctahedron obtained by permuting the entries and considering all sign flips of the point $(\pm1/\sqrt{2},\pm1/\sqrt{2},0)$, which form the $\mathbb{O}$-orbit of the vertex $V_2\in\mathcal{R}$. These are stable equilibria of the regularised reduced system \eqref{regredsys} which correspond to collision configurations of the full system \eqref{system} where twelve pairs of vortices simultaneously collide at the vertices of an $\mathbb{S}^2$-inscribed regular cuboctahedron (see Figure \ref{fig: oct collision eqs (b)}).

\item[(c)] The \textbf{cubic collision equilibria} occurring at the vertices of a cube $(\pm1/\sqrt{3},\pm1/\sqrt{3},\allowbreak\pm1/\sqrt{3})$, which form the $\mathbb{O}$-orbit of the vertex $V_3\in\mathcal{R}$. These are stable equilibria of the regularised reduced system \eqref{regredsys} which correspond to collision configurations of the full system \eqref{system} where eight triplets of vortices simultaneously collide at the vertices of an $\mathbb{S}^2$-inscribed regular cube (see Figure \ref{fig: oct collision eqs (c)}).

\end{itemize}
\end{itemize}
\end{thm}

\begin{figure}\centering
\begin{overpic}[width=.5\textwidth]{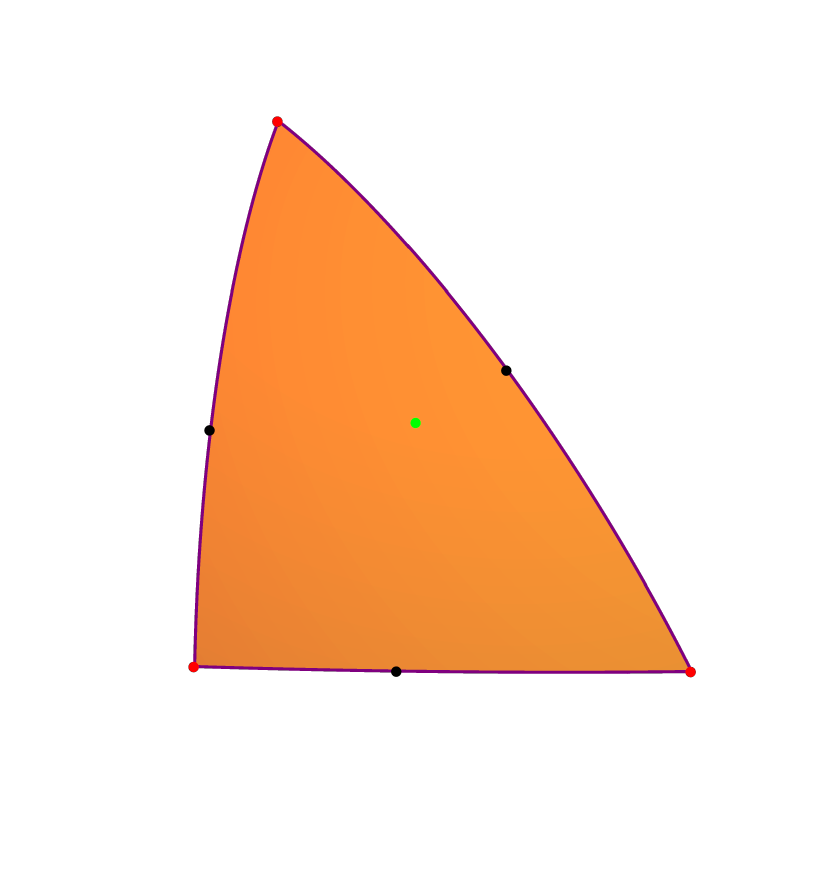}
    \put(17,21){\fontsize{7}{8.4}\selectfont\textcolor{red}{$V_2$}}
    \put(80,21){\fontsize{7}{8.4}\selectfont\textcolor{red}{$V_3$}}
    \put(30,89){\fontsize{7}{8.4}\selectfont\textcolor{red}{$V_1$}}
    
    \put(0,50){\parbox{20mm}{\fontsize{7}{8.4}\selectfont\centering Irreg. Trunc.\newline Octahedron}}
    
    \put(57,57){\parbox{20mm}{\fontsize{7}{8.4}\selectfont\centering Pseudo\newline Rhombicub.}}

    \put(35,45){\parbox{20mm}{\fontsize{7}{8.4}\selectfont\centering \textcolor{green}{Pseudo Snub}\newline \textcolor{green}{Cube (R)}}}
    
    \put(34,17){\parbox{20mm}{\fontsize{7}{8.4}\selectfont\centering Irreg. Trunc.\newline Cube}}

    \end{overpic}
  \caption{Equilibria of the regularised reduced system \eqref{regredsys} on the fundamental region $\mathcal{R}$ (cf. Theorem \ref{thm: Oct eqs}).}
    \label{fig: FullOctregwitheqs}
\end{figure}

A detailed discussion of the edge lengths and truncation ratios of the polyhedra in Theorem \ref{thm: Oct eqs} is given in Appendix A. The equilibria of the regularised reduced system \eqref{regredsys} on the fundamental region $\mathcal{R}$ are illustrated in Figure \ref{fig: FullOctregwitheqs}.

\begin{rem}\label{rem: stab/instab in full system} A consequence of the theorem is that the instability of the equilibria in item (i) (b), (c), (d) of Theorem \ref{thm: Oct eqs} in the regularised reduced system \eqref{regredsys} translates into the instability of the corresponding configurations in the full system \eqref{system}, depicted in Figure \ref{fig: oct non collision eqs (b)}, \ref{fig: oct non collision eqs (c)}, \ref{fig: oct non collision eqs (d)}. It would be interesting to investigate the stability of the configurations in the full system \eqref{system} corresponding to the pseudo snub cube equilibria in item (i) (a) (see Figure \ref{fig: oct non collision eqs (a)}).
\end{rem}

\begin{proof}[Proof of Theorem \ref{thm: Oct eqs}:] It was already observed that the collision equilibria are stable and occur at the points in $\mathcal{F}[\mathbb{O}]$. These points coincide with the union of the three distinct $\mathbb{O}$-orbits of $V_1$, $V_2$, $V_3$. Then, item (ii) of the theorem is an immediate consequence of these observations and  Proposition \ref{prop: n-uple collisions}. As for item (i) of the theorem, in this proof we only focus on proving the existence of  the equilibria of the regularised reduced system \eqref{regredsys} described in items (a), (b), (c), (d): their uniqueness can be proved with elementary methods, through direct estimates, by studying the monotonicity of the polynomials involved (the reader might refer to Lemma 3.3.3 and Lemma 3.3.4 of \cite{Poma}).
We begin with the edges (excluding the vertices): from the discussion in section \ref{subsec: edges fund reg}, we know that we should search for the critical points of the restrictions $\restr{p_{\mathbb{O}}}{\mathcal{E}_{ij}}$. 
    \begin{itemize}
    
    \item Edge $\mathcal{E}_{12}$. Ruling out $V_1$ and $V_2$, this edge can be parametrised as follows: \[\mathcal{E}_{12}=\{(x,0,\sqrt{1-x^2})\colon 0<x<1/\sqrt{2}\}\,,\] and a direct computation shows that:
    \begin{align*}
    &\frac{d}{dx}p_{\mathbb{O}}(x,0,\sqrt{1-x^2})=f_{1,2}(x)\cdot M_{1,2}(x)\,\text{  with:}\\
   &f_{1,2}(x)\coloneqq-6 + 
   41 x^2 - 72 x^4 + 22 x^6 + 89 x^8 - 120 x^{10} + 40 x^{12}\,.
   \end{align*}
   The complete expression of the polynomial $M_{1,2}(x)$ can be easily read from \eqref{octedgeqpoly1}.  The factors appearing in $M_{1,2}(x)$ do not vanish in the interval $(0,\,1/\sqrt{2})$ and $f_{1,2}$ has exactly one root $\bar{x}$ in the same interval (see the proof of Lemma 3.3.4 in \cite{Poma} for technical details). A numerical approximation of such a root is $\bar{x}\approx0.474726$, which yields $\beta=\bar{x}^2\approx0.225365$.
   
   \item Edge $\mathcal{E}_{13}$. Ruling out $V_1$ and $V_3$, this edge can be parametrised as follows: \[\mathcal{E}_{13}=\{(x,x,\sqrt{1-2x^2})\,\colon\, 0<x<1/\sqrt{3}\}\,,\] and a direct computation shows that:
     \begin{align*}
   &\frac{d}{dx}p_{\mathbb{O}}(x,x,\sqrt{1-2x^2})=f_{1,3}(x)\cdot M_{1,3}(x)\, \text{     with:}\\
   &f_{1,3}(x)\coloneqq-3 + 41 x^2 - 
   190 x^4 + 502 x^6 - 323 x^8 - 1023 x^{10} + 1188 x^{12} 
   \end{align*}
   The complete expression of the polynomial $M_{1,3}(x)$ can be easily read from \eqref{octedgeqpoly2}. The factors appearing in $M_{1,3}(x)$ do not vanish in the interval $(0,\,1/\sqrt{3})$ and $f_{1,3}$ has exactly one root $\tilde{x}$ in the same interval (see the proof of Lemma 3.3.4 in \cite{Poma} for technical details). A numerical approximation of such a root is $\tilde{x}\approx0.350294$, which yields $\gamma=\tilde{x}^2\approx0.122706$.
   
   \item Edge $\mathcal{E}_{23}$. Ruling out $V_2$ and $V_3$, this edge can be parametrised as follows: \[\mathcal{E}_{23}=\{(x,\sqrt{1-2x^2},x)\,\colon\, 1/\sqrt{3}<x<1/\sqrt{2}\}\,,\] and a direct computation shows that:
   \[ \frac{d}{dx}p_{\mathbb{O}}(x,\sqrt{1-2x^2},x)=f_{1,3}(x)\cdot M_{1,3}(x) \]
   The complete expression of the polynomial $M_{1,3}(x)$ can be easily read from \eqref{octedgeqpoly2}, as for the edge $\mathcal{E}_{13}$. The factors appearing in $M_{1,3}(x)$ do not vanish in the interval $(1/\sqrt{3},\,1/\sqrt{2})$ and $f_{1,3}$ has exactly one root $\hat{x}$ in the same interval (see the proof of Lemma 3.3.4 in \cite{Poma} for technical details). A numerical approximation of such a root is $\hat{x}\approx0.685580$, which yields $\delta=\hat{x}^2\approx0.470019$.
    \end{itemize}
    
This proves items (b), (c), (d) of item (i), except for the instability property.
Next, we continue with the interior $\mathrm{int}(\mathcal{R})$. In view of \eqref{reg red oct ham in the inv}, the critical points in the interior $\mathrm{int}(\mathcal{R})$ are the real solutions $(s,t)\in U$ of the system:
\begin{equation*}
\begin{split}
\frac{\partial}{\partial s}p^I_{\mathbb{O}}(s,t)=
F_1(s,t)\, F_2(s,t)^2G_1(s,t)=0\\
\frac{\partial}{\partial t}p^I_{\mathbb{O}}(s,t)=
F_1(s,t)\, F_2(s,t)^2H_1(s,t)=0
\end{split}
\end{equation*}
where the complete expressions of the polynomials $G_1$, $H_1$ are given in \eqref{octdsfactor} and \eqref{octdtfactor}. 
The $0$ level set of $p_{\mathbb{O}}$ coincides with the $0$ level set of $\tilde{h}_{\mathbb{O}}$, which contains only collision equilibria. There is no collision equilibrium in $\mathrm{int}(\mathcal{R})$, therefore $F_1$, $F_2$ do not simultaneously vanish in $U$ and can be ignored. We are left with two polynomials $G_1(s,t)$, $H_1(s,t)$ which might contain common real zeroes in $U$. Therefore, the problem of finding critical points of $p_{\mathbb{O}}$ in $\mathrm{int}(\mathcal{R})$ is reduced to solving $G_1(s,t)=0$, $H_1(s,t)=0$ for $(s,t)\in U\subset(0,\,1/3)\times(0,\,1/27)$ (see Remark \ref{rem: oct inv fund reg bounds}). The $s$- and $t$-values of these common zeroes of $G_1$, $H_1$ are real roots in these intervals of the resultants $\mathrm{Res}_t(G_1, H_1)$ and $\mathrm{Res}_s(G_1, H_1)$, respectively. These resultants allow a factorisation of the form:
\[\mathrm{Res}_s(G_1,H_1)(t)=q_2(t)\cdot L_2(t)\,,\]
\[\mathrm{Res}_t(G_1,H_1)(s)=q_1(s)\cdot L_1(s)\,.\] 
The complete expressions of these resultants and of the polynomials $q_1(s)$ and $q_2(t)$ are given in \eqref{octinteqpoly2}, \eqref{octresultants} and  \eqref{octinteqpoly1},  \eqref{octresultantt}. The expressions for the polynomials $L_1(s)$ and $L_2(
t)$ can easily be read from \eqref{octresultants} and \eqref{octresultantt}. The zeroes of the factors $L_2(t)$ and $L_1(s)$ do not lead to any common root $(s,t)\in U$ of $G_1$, $H_1$, so we focus on the polynomials $q_1(s)$ and $q_2(t)$: $q_1$ has a unique root $\bar{s}$ satisfying $0<\bar{s}<1/3$, while $q_2$ has a unique root $\bar{t}$ satisfying $0<\bar{t}<1/27$ (see the proof of Lemma $3.3.3$ in \cite{Poma} for technical details). Numerical evidence\footnote{The restriction  $\restr{p_{\mathbb{O}}}{\mathcal{R}}$ to the compact domain $\mathcal{R}$ must attain its global maximum at some critical point, and it can be verified that $\max_{p\in\partial\mathcal{R}}p_{\mathbb{O}}<p^I_{\mathbb{O}}(\tilde{s},\tilde{t})$ where $(\tilde{s},\tilde{t})$ is a numerical approximation of $(\bar{s},\bar{t})$.} confirms that $(\bar{s},\bar{t})$ is a common zero of $H_1$ and $G_1$ in $U$ (see Figure \ref{fig: FullOctInvReg}). The point $(\bar{s},\bar{t})$ in the invariant plane is approximately $(0.202686,\,0.009042)$. 

\begin{figure}
   \centering
   \begin{overpic}[width=0.98\textwidth]{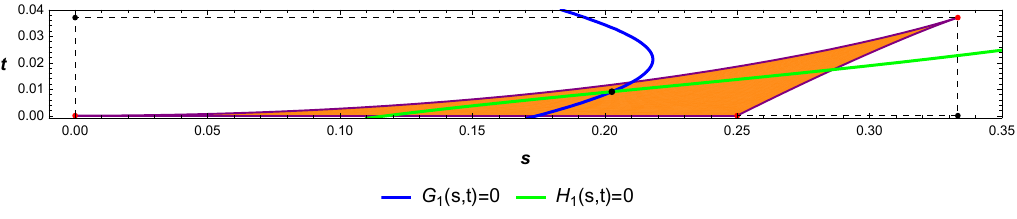}
    \put(8,11){\textcolor{red}{$V_1$}}
    \put(69,11){\textcolor{red}{$V_2$}}
    \put(95,18){\textcolor{red}{$V_3$}}
    \put(88,7){\textcolor{black}{$(1/3,\,0)$}}
    \put(8,17){\textcolor{black}{$(0,\,1/27)$}}
    \put(55,14){\textcolor{black}{$(\bar{s},\bar{t})$}}
    \end{overpic}
    \caption{In orange, the domain $U=\psi^{-1}(\mathrm{int}(\mathcal{R}))$ in the invariants $(s, t)$-plane that parametrises the interior $\mathrm{int}(\mathcal{R})$ of the fundamental region. In blue and green the local vanishing loci of $G_1$, $H_1$.}
    \label{fig: FullOctInvReg}
\end{figure}

Once we have found the equilibria, in order to determine their (Lyapunov) stability/instability it is enough to parametrise the northern hemisphere (that contains the fundamental region $\mathcal{R}$) as $z=\sqrt{1-x^2-y^2}$ and evaluate the Hessian of $P(x, y)\coloneqq p_{\mathbb{O}}(x, y, \sqrt{1-x^2-y^2})$ at the critical points in  $\mathcal{R}$ (excluding $V_1$, $V_2$ and $V_3$); the remaining critical points have the same stability/instability properties by $\mathbb{O}_h$-symmetry. We exhibit a numerical evidence of the stability properties, by evaluating the Hessian at numerical approximations of the found non-collision equilibria\footnote{In fact, if the Hessian is positive/negative definite or indefinite at a point $p\in\R^2$, it will be such at any point $p'$ in an open neighbourhood of $p$.}:
\begin{equation*}\begin{aligned} &\mathrm{Hess}(P)(\alpha_1,\alpha_2)\approx\begin{pmatrix}-0.151878&-0.092110\\-0.092110&-0.094324\end{pmatrix}\quad\text{Negative definite}\,,\\&\mathrm{Hess}(P)(\sqrt{\beta},0)\approx\begin{pmatrix}-0.268813&0\\0&0.046149\end{pmatrix}\quad\text{Indefinite}\,,\\&\mathrm{Hess}(P)(\sqrt{\gamma},\sqrt{\gamma})\approx\begin{pmatrix}-0.126697&-0.140807\\-0.140807&-0.126697\end{pmatrix}\quad\text{Indefinite}\,,\\&\mathrm{Hess}(P)(\sqrt{\delta},\sqrt{1-2\delta})\approx\begin{pmatrix}0.151531&0.027061\\0.027061&-0.195431\end{pmatrix}\quad\text{Indefinite}\,.\end{aligned}
\end{equation*}
\end{proof}

\subsection{Reduced phase portrait}\label{sec: oct dynamics, ph portrait}
\begin{figure}
    \centering
    \renewcommand{\figurename}{Figure}\includegraphics[width=.85\linewidth]{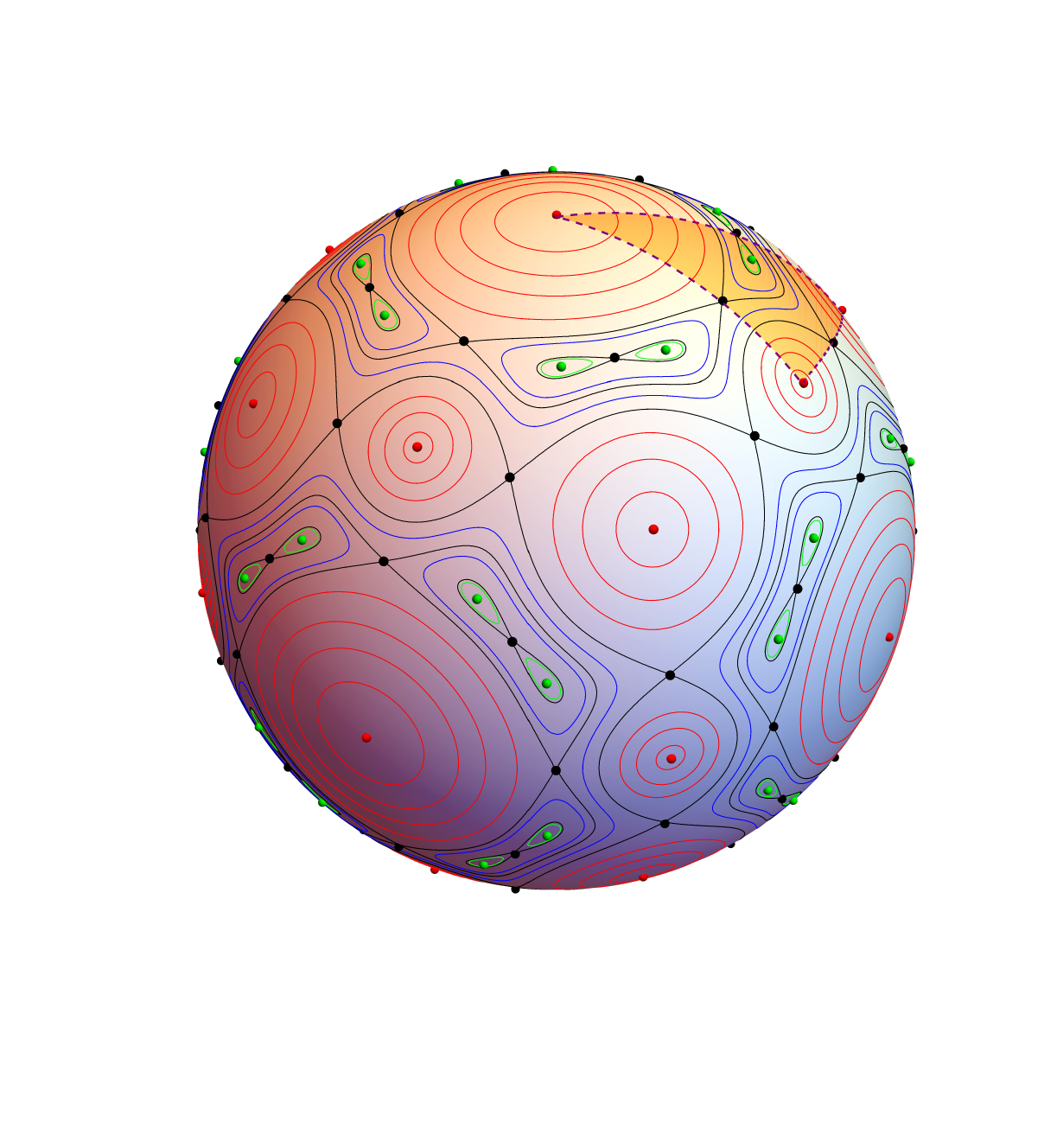}
    \caption{The (regularised) reduced phase portrait of octahedrally symmetric solutions of the $24$-vortex problem on the sphere, with the fundamental region $\mathcal{R}$ highlighted.}
    \label{fig: OctPhasePortrait}
\end{figure}

Finally, we present in Figure \ref{fig: OctPhasePortrait} a numerical approximation of the phase portrait of the regularised reduced system \eqref{regredsys}. The pseudo snub cube equilibrium points are indicated in green, the collision equilibrium points (cube, octahedron, cuboctahedron) are indicated  in red and the pseudo rhombicuboctahedron, truncated octahedron and truncated cube equilibrium points are indicated in black. The colour code is for periodic orbits near the stable equilibria or heteroclinic/homoclinic orbits emanating from the unstable equilibria. Blue periodic orbits are other families of periodic orbits not emanating from equilibria.
\section{Icosahedrally symmetric solutions of the $60$-vortex problem on the sphere}

Finally, we turn to the case $K=\mathbb{I}$ and $N=\abs{K}=60$. Again, we consider a slight variation of the regularised Hamiltonian $\tilde{h}_{\mathbb{I}}$:
\[p_{\mathbb{I}}(x,y,z)\coloneqq 2^{34}\sqrt[30]{\tilde{h}_{\mathbb{I}}(x,y,z)}\,,\]
which has the same critical points and the same level sets as $\tilde{h}_{\mathbb{I}}$ (corresponding to different values, though). The complete expression of $\tilde{h}_{\mathbb{I}}$ is given by \eqref{icoregham}, while the complete expression of $p_{\mathbb{I}}$ can be read from \eqref{icosimplregham}.

The full icosahedral group $\mathbb{I}_h$ is generated by the following reflections:
\begin{align*}
    &r_{e_x}=\begin{pmatrix}
-1 & 0 &0\\
0 & 1&0\\
0&0&1
\end{pmatrix}\,,\\
&r_{e_y}=\begin{pmatrix}
1 & 0 &0\\
0 & -1&0\\
0&0&1
\end{pmatrix}\,,\\
&r_{-\phi^2\,e_x-\phi\, e_y+e_z}=\frac{1}{2}\begin{pmatrix}
1-\phi & -\phi &1\\
-\phi & 1&\phi-1\\
1&\phi-1&\phi
\end{pmatrix}\,,\\&\text{where: }\phi\coloneqq \frac{1+\sqrt{5}}{2}\,.
\end{align*}
The reflection planes of such generators contain the three edges of the fundamental region $\mathcal{R}$ in green in Figure \ref{fig: Full Ico Fund Reg}. Explicitly, $\mathcal{R}$ is described by:
\[\mathcal{R}\coloneqq\{(x,y,z)\in\mathbb{S}^2\,\colon\,y\ge 0\,\wedge\,x\ge 0 \,\wedge\,z\ge\phi^2\,x+\phi\,y \}\,.\]
We denote by $\mathcal{E}_{12}$ the edge $\mathcal{R}\cap H_{e_y}$, by $\mathcal{E}_{13}$ the edge $\mathcal{R}\cap H_{e_x}$ and by  $\mathcal{E}_{23}$ the edge $\mathcal{R}\cap H_{-\phi^2\,e_x-\phi\,e_y+e_z}$. The three vertices of $\mathcal{R}$ are:
\[\begin{aligned} &V_1\coloneqq \mathcal{E}_{12}\cap\mathcal{E}_{13}=(0,0,1)\,,\\& V_2\coloneqq \mathcal{E}_{12}\cap\mathcal{E}_{23}=\frac{1}{\sqrt{9\phi+6}}(\phi,\,0,\,1+2\phi)\,,\\&V_3\coloneqq\mathcal{E}_{13}\cap\mathcal{E}_{23}=\frac{1}{\sqrt{2+\phi}}(0,\,1,\,\phi)\,.\end{aligned}\]

Finally, we describe the isotropy group of points on $\partial\mathcal{R}$: these are the only points with non-trivial isotropy group. The isotropy group of any point $r\in\mathcal{E}_{12}\setminus\{V_1,V_2\}$ is $(\mathbb{I}_h)_r=\langle \,r_{e_y}\,\rangle$; similarly, for any point $r\in\mathcal{E}_{13}\setminus\{V_1,V_3\}$ we have $(\mathbb{I}_h)_r=\langle \,r_{e_x}\,\rangle$ and for any point $r\in\mathcal{E}_{23}\setminus\{V_2,V_3\}$ we have $(\mathbb{I}_h)_r=\langle\, r_{-\phi^2\,e_x-\phi\,e_y+e_z}\,\rangle$. As for the vertices, $(\mathbb{I}_h)_{V_1}=\langle\, r_{e_y},\,r_{e_x}\,\rangle$, $(\mathbb{I}_h)_{V_2}=\langle\, r_{e_y},\,r_{-\phi^2\,e_x-\phi\,e_y+e_z}\,\rangle$ and $(\mathbb{I}_h)_{V_3}=\langle\, r_{e_x},\,r_{-\phi^2\,e_x-\phi\,e_y+e_z}\,\rangle$.

\subsection{The polynomial $p_{\mathbb{I}}$ in invariant coordinates}

Before dealing with the search for equilibria, we rename the icosahedral fundamental invariants of Theorem \ref{thm: CST} (their complete expressions are \eqref{icoinv6}, \eqref{icoinv10}):  \begin{align*} &s\coloneqq \theta_2(x,y,z)\,,\\ &t\coloneqq \theta_3(x,y,z)\,.\end{align*}
Using the parametrisation $\psi$ in \eqref{param}, we obtain the polynomial function $p^I_{\mathbb{I}}\colon U\to\R$ in the variables $(s,t)$.
An explicit expression of $p^I_{\mathbb{I}}$ is the following (see Subsection 4.2.1 of \cite{Poma} for details on its derivation):
\begin{equation}\label{reg red ico ham in the inv}\begin{aligned}&p^I_{\mathbb{I}}(s,t)=c\,F_4(s,t)^4 F_5(s,t)^2 F_6(s,t),\,\text{where:}\\&c\coloneqq-2^{122}\cdot3^{-14}\cdot5^{-12}\\
   &F_4(s,t)\coloneqq -1 - 12 s + 45 s^2 + 27 t\,,\\
   &F_5(s,t)\coloneqq 649 - 7524 s + 14256 s^2 + 13500 s^3 + 1674 t - 4212 s t + 81 t^2\,,\\
   &\begin{aligned}F_6(s,t)\coloneqq&+5729 - 103566 s + 569988 s^2 - 1118178 s^3 + 1624050 s^4 + 
 24300 s^5 \\ &+ 19701 t - 142236 s t - 32076 s^2 t - 85050 s^3 t + 
 14823 t^2 + 24786 s t^2 + 243 t^3\,.\end{aligned}
   \end{aligned}\end{equation}

\begin{figure}
   \centering
   \begin{overpic}[width=0.98\textwidth]{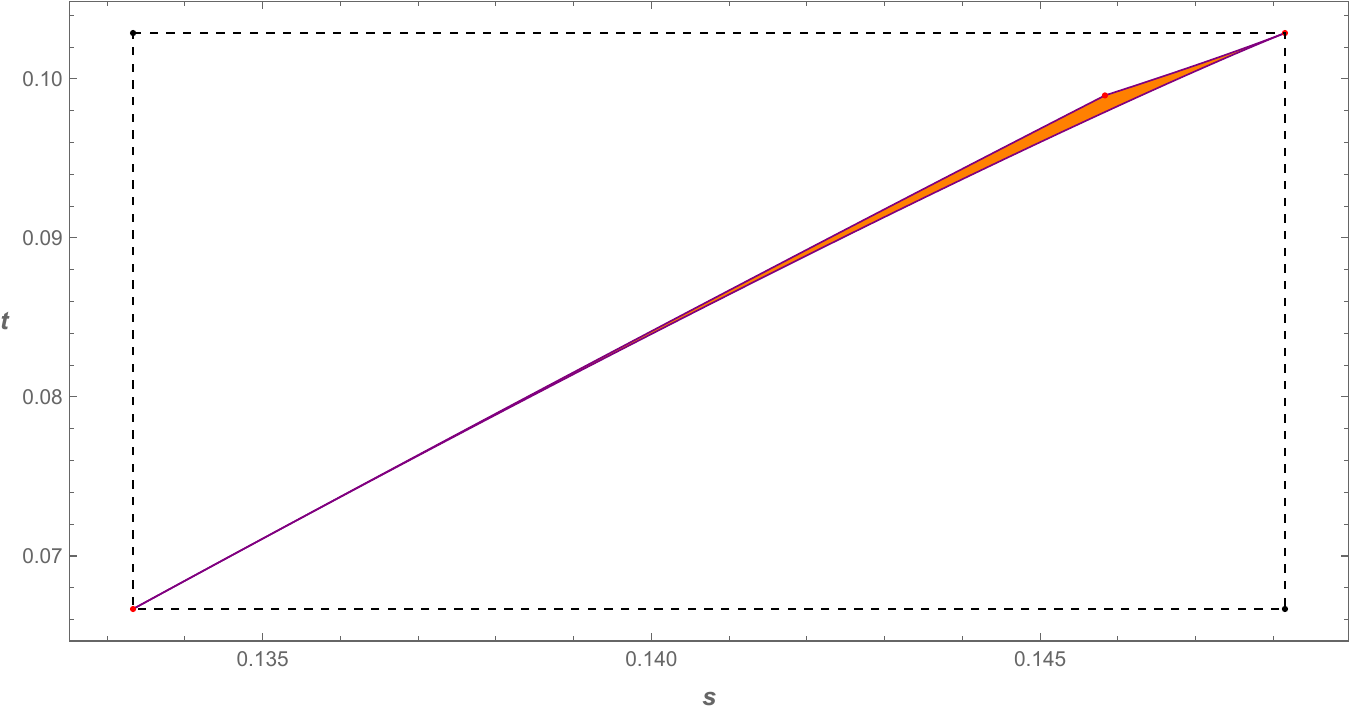}
    \put(6,7){\textcolor{red}{$V_
    3$}}
    \put(80,47){\textcolor{red}{$V_1$}}
    \put(96,50){\textcolor{red}{$V_2$}}
    \put(80, 9){\textcolor{black}{$(4/27,\,1/15)$}}
    \put(11,48){\textcolor{black}{$(2/15,\,25/243)$}}
    \end{overpic}
    \caption{In orange, the domain $U=\psi^{-1}(\mathrm{int}(\mathcal{R}))$ in the invariants $(s, t)$-plane that parametrises the interior $\mathrm{int}(\mathcal{R})$ of the fundamental region.}
    \label{fig: FullIcoInvRegNoGH}
\end{figure}

\begin{rem}\label{rem: ico inv fund reg bounds} Let $(x,y,z)$ be any point in $\mathrm{int}(\mathcal{R})$ and $(s,t)=\psi^{-1}(x,y,z)$ be the corresponding point in $U$. Unlike in the octahedral case, this time the algebraic complexity of the invariants hinders any attempt to find precise bounds on $s$ and $t$ through inequalities. Nevertheless, Figure \ref{fig: FullIcoInvRegNoGH} suggests taking \[\begin{aligned} &\min_{i=1,2,3}\{\theta_2(V_i)\}<s<\max_{i=1,2,3}\{\theta_2(V_i)\}\,,\\& \min_{i=1,2,3}\{\theta_3(V_i)\}<t<\max_{i=1,2,3}\{\theta_3(V_i)\}\,,\end{aligned}\] yielding $2/15<s<4/27$ and $1/15<t<25/243$. This shows heuristically that $U\subseteq (2/15,\,4/27)\times (1/15,\,25/243)$. Even though we did not prove these estimates, we are going to use them in the proof of Theorem \ref{thm: Ico eqs}, since in any case we will resort to numerical tools.
\end{rem}

\subsection{Icosahedrally symmetric collision and equilibrium configurations}\label{subsec: ico eqs}
\begin{figure}
    \centering
    \subfloat[Pseudo snub dodecahedron equilibrium configuration.\label{fig: ico non collision eqs (a)}]{\includegraphics[width=0.4\linewidth]{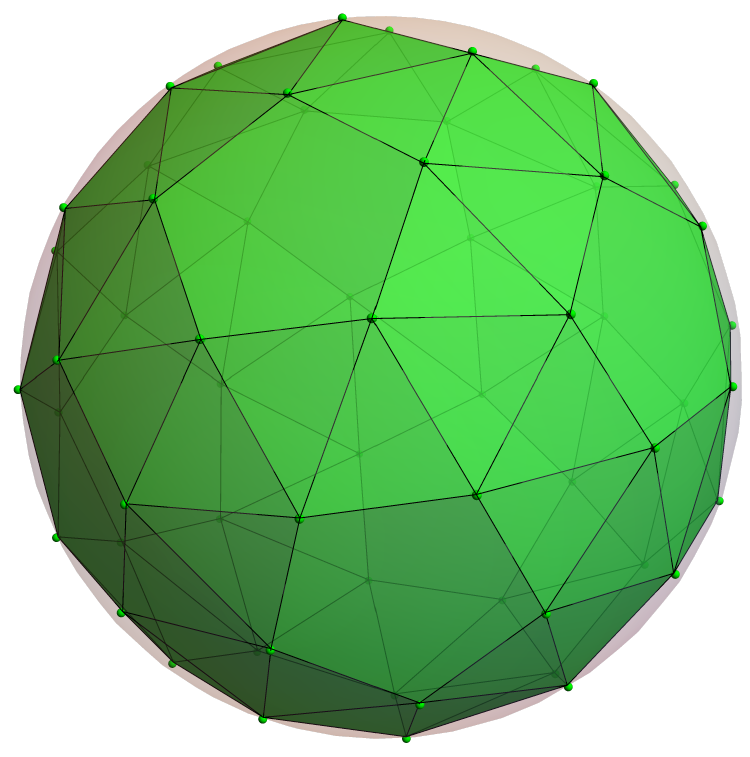}}
        \quad
    \subfloat[Truncated dodecahedron equilibrium configuration.\label{fig: ico non collision eqs (b)}]{\includegraphics[width=0.4\linewidth]{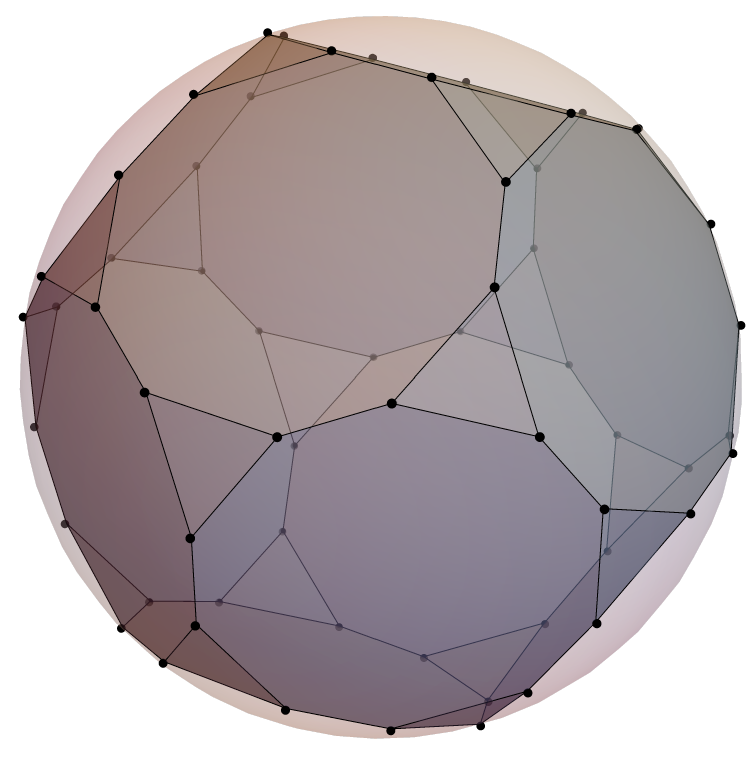}}
    \quad
    \subfloat[Truncated icosahedron equilibrium configuration.\label{fig: ico non collision eqs (c)}]{\includegraphics[width=0.4\linewidth]{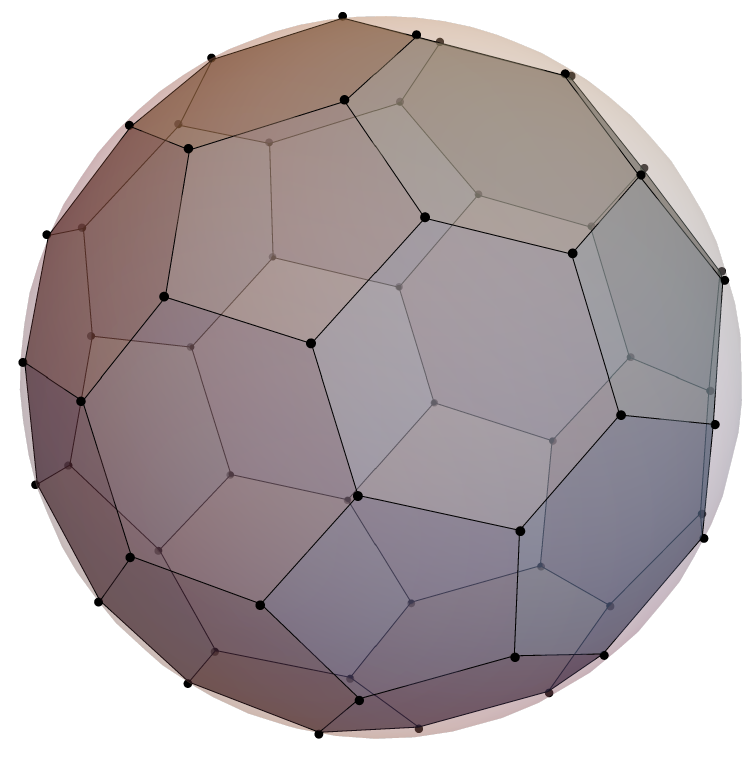}}
    \quad
    \subfloat[Pseudo rhombicosidodecahedron equilibrium configuration.\label{fig: ico non collision eqs (d)}]{\includegraphics[width=0.4\linewidth]{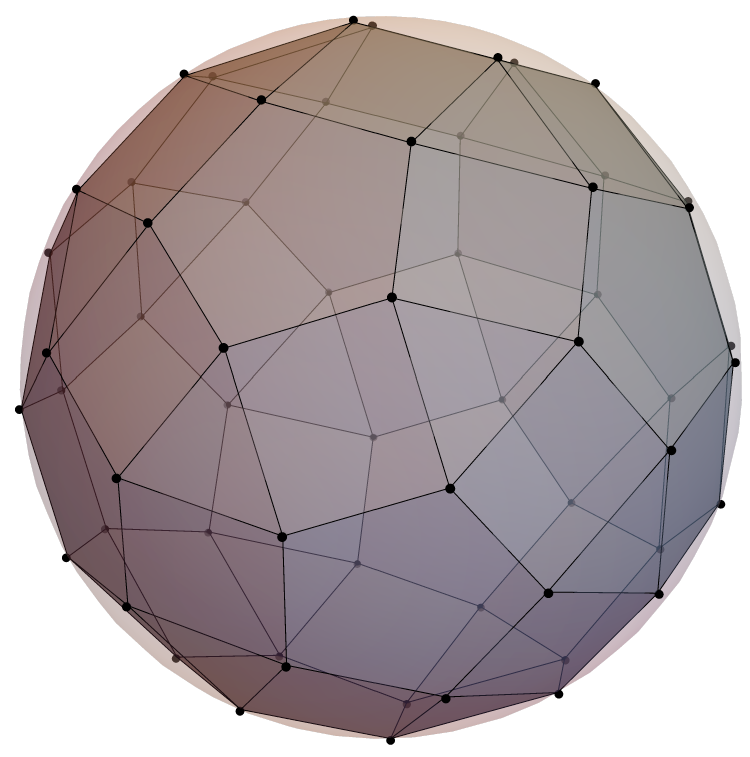}}
    \caption{Equilibrium configurations of the full system \eqref{system} described in item (i) of Theorem \ref{thm: Ico eqs}. For more details on the dimensions, see Appendix A.}\label{fig: ico non collision eqs}
\end{figure}

\begin{figure}
    \centering
    \subfloat[Icosidodecahedral collision configuration: simultaneous binary collision at each vertex.\label{fig: ico collision eqs (a)}]{\includegraphics[width=0.4\linewidth]{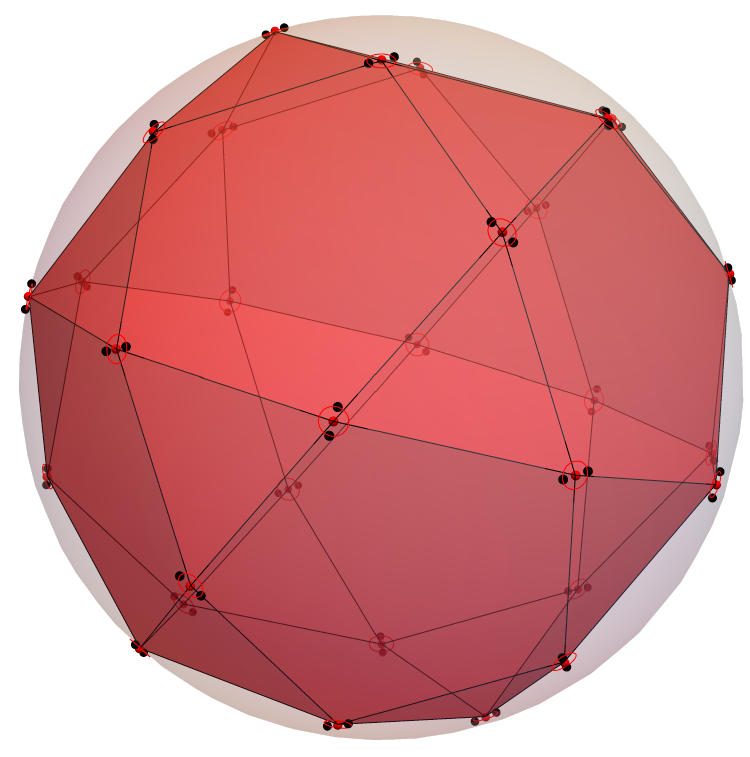}}
        \quad
    \subfloat[Dodecahedral collision configuration: simultaneous ternary collision at each vertex.\label{fig: ico collision eqs (b)}]{\includegraphics[width=0.4\linewidth]{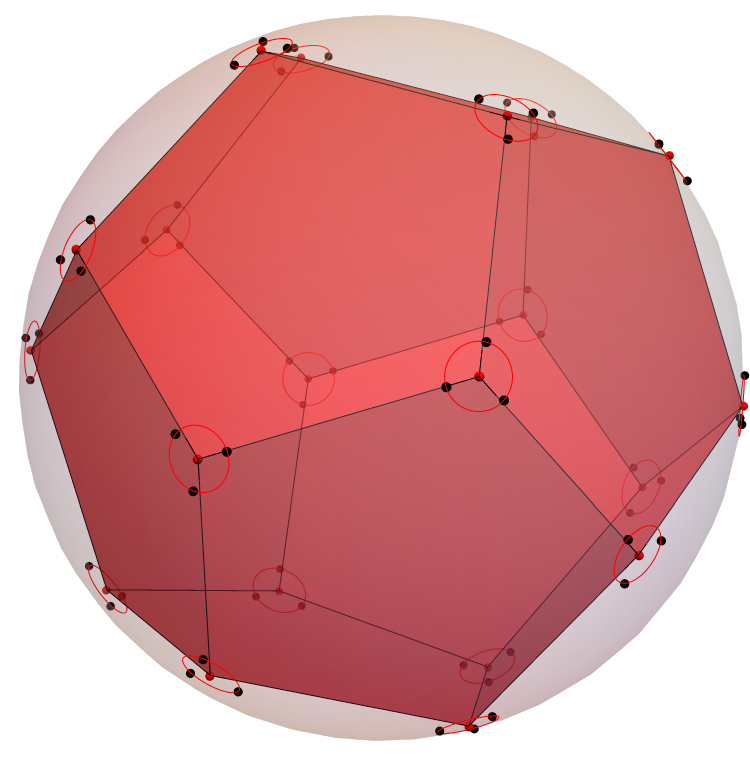}}
    \quad
    \subfloat[Icosahedral collision configuration: simultaneous quinary collision at each vertex.\label{fig: ico collision eqs (c)}]{\includegraphics[width=0.4\linewidth]{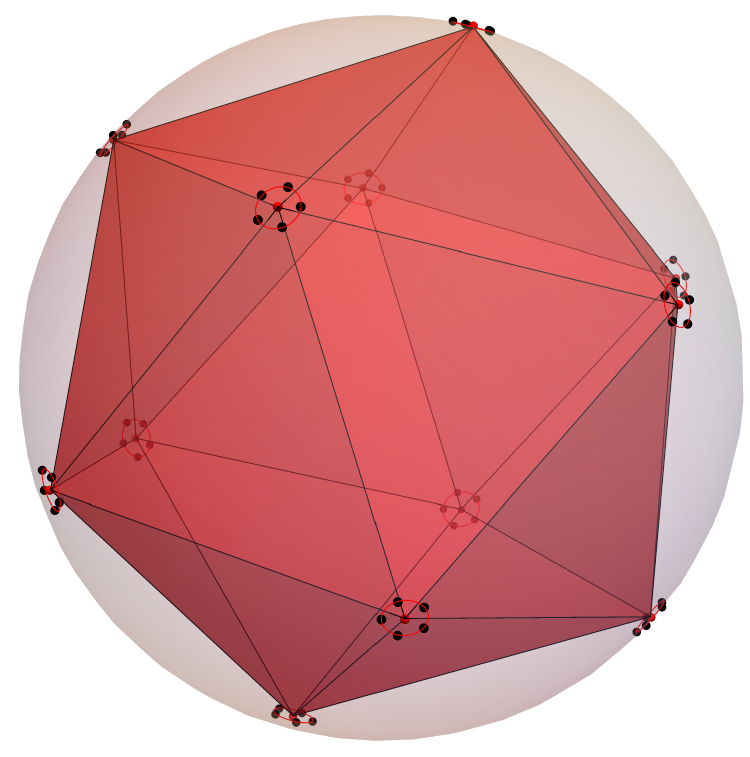}}
    \caption{Collision configurations of the full system \eqref{system} described in item (ii) of Theorem \ref{thm: Ico eqs}.} \label{fig: ico collision eqs}
\end{figure}

We present now the main theorem of this section:
\begin{thm}[Icosahedrally symmetric equilibrium and collision configurations] \label{thm: Ico eqs}  Let $K=\mathbb{I}$, and $N=60$. 
\begin{itemize}
\item[(i)]The following are non-collision equilibria of the regularised reduced system \eqref{regredsys}:
\begin{itemize}

    \item[(a)] The \textbf{pseudo}\footnote{We use the word “pseudo'' to point out that it has the same vertex configuration but irregular faces.} \textbf{snub dodecahedron equilibria} occurring at the two distinct $\mathbb{I}$-orbits contained in the $\mathbb{I}_h$-orbit of the point: \[(\xi_1,\xi_2,\xi_3)\in \mathrm{int}(\mathcal{R})\]
    whose coordinates satisfy: \[\theta_2(\xi_1,\xi_2,\xi_3)=\bar{s},\quad \theta_3(\xi_1,\xi_2,\xi_3)=\bar{t}\]
     where $\bar{s}$ is a root in the interval $(2/15,\,4/27)$ of a polynomial $q_3(s)$ of degree $23$ (see \eqref{icointeqpoly1} for its complete expression) and $\bar{t}$ is a root in the interval $(1/15,\,25/243)$ of a polynomial $q_4(t)$ of degree $23$ (see \eqref{icointeqpoly2} for its complete expression). These are stable equilibria of the regularised reduced system \eqref{regredsys} which correspond to equilibrium configurations of the full system \eqref{system} where the vortices are situated at the vertices of an $\mathbb{S}^2$-inscribed left-chiral pseudo snub dodecahedron or its right-chiral version (see Figure \ref{fig: ico non collision eqs (a)} for the right-chiral version). The point $(\xi_1,\xi_2,\xi_3)$ defined above is approximately given by $(0.155679, 0.206697, 0.965940)$, the root $\bar{s}$ is approximately given by $0.143135$ and the root $\bar{t}$ is approximately given by $0.091740$.
    
    \item[(b)] The \textbf{truncated dodecahedron equilibria} occurring at the $\mathbb{I}$-orbit of the point:
    \[(\sqrt{\zeta},\,0,\,\sqrt{1-\zeta})\in \mathcal{E}_{12}\]
    where $\zeta$ is a certain root in the interval $(0,\,\phi^2/(6+9\phi))$ of a polynomial $g_{1,2}(u)$ of degree $15$ (see \eqref{Icopolyedge12} for its complete expression). These are unstable equilibria of the regularised reduced system \eqref{regredsys} which correspond to equilibrium configurations of the full system \eqref{system} where the vortices are situated at the vertices of an $\mathbb{S}^2$-inscribed irregular truncated dodecahedron (see Figure \ref{fig: ico non collision eqs (b)}). The point $(\sqrt{\zeta},\,0,\,\sqrt{1-\zeta})$ defined above is approximately given by $(0.144420, 0, 0.989517)$.

    \item[(c)] The \textbf{truncated icosahedron equilibria} occurring at the $\mathbb{I}$-orbit of the point:
    \[(0,\,\sqrt{\eta},\,\sqrt{1-\eta})\in \mathcal{E}_{13}\]
    where $\eta$ is a certain root in the interval $(0,\,1/(2+\phi))$ of a polynomial $g_{1,3}(v)$ of degree $15$ (see \eqref{Icopolyedge13} for its complete expression). These are unstable equilibria of the regularised reduced system \eqref{regredsys} which correspond to equilibrium configurations of the full system \eqref{system} where the vortices are situated at the vertices of an $\mathbb{S}^2$-inscribed irregular truncated icosahedron (see Figure \ref{fig: ico non collision eqs (c)}). The point $(0,\,\sqrt{\eta},\,\sqrt{1-\eta})$ defined above is approximately given by $(0, 0.191446, 0.981503)$.

    \item[(d)] The \textbf{pseudo rhombicosidodecahedron equilibria} occurring at the $\mathbb{I}$-orbit of the point:
    \[\frac{(1-\bar{\lambda}\,)(\phi,\,0,\,1+2\phi)+\bar{\lambda}(0,\,1,\,\phi)}{\norm{(1-\bar{\lambda}\,)(\phi,\,0,\,1+2\phi)+\bar{\lambda}(0,\,1,\,\phi)}}\coloneqq k(\bar{\lambda})\in \mathcal{E}_{23}\]
    where $\bar{\lambda}$ is a certain root in the interval $(0,\,1)$ of a polynomial $g_{2,3}(\lambda)$ of degree $30$ (see \eqref{Icopolyedge23} for its complete expression). These are unstable equilibria of the regularised reduced system \eqref{regredsys} which correspond to equilibrium configurations of the full system \eqref{system} where the vortices are situated at the vertices of an $\mathbb{S}^2$-inscribed pseudo rhombicosidodecahedron (see Figure \ref{fig: ico non collision eqs (d)}). The point $k(\bar{\lambda})$ defined above is approximately given by $(0.213451, 0.239928, 0.947034)$.
\end{itemize}

\item[(ii)]The only collision equilibria of the regularised reduced system \eqref{regredsys} are:

\begin{itemize}

\item[(a)] The \textbf{icosidodecahedral collision equilibria} occurring at the vertices of an icosidodecahedron (obtained considering all the even permutations and all sign flips of the points $(\pm 1,\,0,\,0)$ and $1/2 (\pm 1,\,\pm(\phi-1),\,\pm\phi) $), which form the $\mathbb{I}$-orbit of the vertex $V_1\in\mathcal{R}$. These are stable equilibria of the regularised reduced system \eqref{regredsys} which correspond to collision configurations of the full system \eqref{system} where thirty pairs of vortices simultaneously collide at the vertices of an $\mathbb{S}^2$-inscribed regular icosidodecahedron (see Figure \ref{fig: ico collision eqs (a)}).
\item[(b)] The \textbf{dodecahedral collision equilibria} occurring at the vertices of a dodecahedron (obtained considering all the even permutations and the sign flips of the point $1/\sqrt{6+9\phi}\,(\pm\phi,\, 0,\,\pm(1+2\phi))$ and all the sign flips of the point $1/\sqrt{3}(\pm 1,\,\pm 1,\,\pm 1)$), which form the $\mathbb{I}$-orbit of the vertex $V_2\in\mathcal{R}$. These are stable equilibria of the regularised reduced system \eqref{regredsys} which correspond to collision configurations of the full system \eqref{system} where twenty triplets of vortices simultaneously collide at the vertices of an $\mathbb{S}^2$-inscribed regular dodecahedron (see Figure \ref{fig: ico collision eqs (b)}).
\item[(c)] The \textbf{icosahedral collision equilibria} occurring at the vertices of an icosahedron (obtained considering all the even permutations and the sign flips of the point  $1/\sqrt{2+\phi}\allowbreak\,(0,\,\pm 1,\,\pm\phi)$), which form the $\mathbb{I}$-orbit of the vertex $V_3\in\mathcal{R}$. These are stable equilibria of the regularised reduced system \eqref{regredsys} which correspond to collision configurations of the full system \eqref{system} where twelve quintuplets of vortices simultaneously collide at the vertices of an $\mathbb{S}^2$-inscribed regular icosahedron (see Figure \ref{fig: ico collision eqs (c)}).

\end{itemize}
\end{itemize}
\end{thm}

As in the previous section, a detailed discussion of the edge lengths and truncation ratios of the polyhedra in Theorem \ref{thm: Ico eqs} is given in Appendix A, while the equilibria are shown in Figure \ref{fig: FullIcoregwitheqs}.

\begin{rem} The analogous conclusions of Remark \ref{rem: stab/instab in full system} about the instability in the full system can be drawn: the configurations in the full system corresponding to equilibria in item (i) (b), (c), (d) of Theorem \ref{thm: Ico eqs} are unstable. On the other hand, it would be interesting to investigate the stability of the configurations in the full system \eqref{system} corresponding to the pseudo snub dodecahedron equilibria in item (i) (a) (see Figure \ref{fig: ico non collision eqs (a)}).
\end{rem}

\begin{proof}[Proof of Theorem \ref{thm: Ico eqs}:] It was already observed that the collision equilibria are stable and occur at the points in $\mathcal{F}[\mathbb{I}]$. These points coincide with the union of the three distinct $\mathbb{I}$-orbits of $V_1$, $V_2$, $V_3$. Then, item (ii) of the theorem is an immediate consequence of these observations and  Proposition \ref{prop: n-uple collisions}.

Now, we study the edges (excluding the vertices). From the discussion in Subsection \ref{subsec: edges fund reg}, we know that we should search for the critical points of the restrictions $\restr{p_{\mathbb{I}}}{\mathcal{E}_{ij}}$. 
    \begin{itemize}
    \item Edge $\mathcal{E}_{12}$. Ruling out $V_1$ and $V_2$, this edge can be parametrised as follows: \[\mathcal{E}_{12}\setminus\{V_1,V_2\}=\{(x,\,0,\,\sqrt{1-x^2})\,\colon\, 0<x<\phi/\sqrt{9\phi+6}\}\,,\] and a direct computation shows that:
\begin{equation*}
    \frac{d}{dx}p_{\mathbb{I}}(x,0,\sqrt{1-x^2})=g_{1,2}(x^2)\cdot N_{1,2}(x)\,.
   \end{equation*}
   For the complete expression of $g_{1,2}$, see \eqref{Icopolyedge12}, while the complete expression of $N_{1,2}$ can be read from \eqref{Icopolyderivedge12}. A numerical analysis detects only one root $\bar{x}$ of $g_{1,2}$ in the interval $(0,\,\phi/\sqrt{9\phi+6})$, and shows that the factors of $N_{1,2}$ do not have roots in this interval. The root $\bar{x}$ is approximately given by $\bar{x}\approx0.144420$, which yields $\zeta=\bar{x}^2\approx0.020857$. 
   
   \item Edge $\mathcal{E}_{13}$. Ruling out $V_1$ and $V_3$, this edge can be parametrised as follows: \[\mathcal{E}_{13}\setminus\{V_1,V_3\}=\{(0,\,y,\,\sqrt{1-y^2})\,\colon\, 0<y<1/\sqrt{2+\phi}\}\,,\] and a direct computation shows that:
     \begin{equation*}
   \frac{d}{dy}p_{\mathbb{I}}(0,y,\sqrt{1-y^2})=g_{1,3}(y^2)\cdot N_{1,3}(y)\,.
   \end{equation*}
   For the complete expression of $g_{1,3}$, see \eqref{Icopolyedge13}, while the complete expression of $N_{1,3}$ can be read from \eqref{Icopolyderivedge13}. A numerical analysis detects only one root $\bar{y}$ of $g_{1,3}$ in the interval $(0,\,1/\sqrt{2+\phi})$, and shows that the factors of $N_{1,3}$ do not have roots in this interval. The root $\bar{y}$ is approximately given by $\bar{y}\approx0.191446$, which yields $\eta=\bar{y}^2\approx0.036652$.
   
   \item Edge $\mathcal{E}_{23}$.
   This edge can be parametrised by normalising the convex combinations of the vertices $V_2$, $V_3$ (geometrically, this corresponds to projecting onto the sphere from the origin the segment $V_2V_3$), giving:
   \begin{align*}&\mathcal{E}_{23}\setminus\{V_2,V_3\}=\left\{k(\lambda)\,\colon\, 0<\lambda<1\right\},\,\,\text{with:}\\&k(\lambda)\coloneqq\frac{(1-\lambda\,)(\phi,\,0,\,1+2\phi)+\lambda(0,\,1,\,\phi)}{\norm{(1-\lambda\,)(\phi,\,0,\,1+2\phi)+\lambda(0,\,1,\,\phi)}}\,,\end{align*} and a direct computation shows that:
   \[ \frac{d}{d\lambda}p_{\mathbb{I}}\left(k(\lambda)\right)=g_{2,3}(\lambda)\cdot N_{2,3}(\lambda) \]
    For the complete expression of $g_{
    2,3}$, see \eqref{Icopolyedge23}, while the complete expression of $N_{2,3}$ can be read from \eqref{Icopolyderivedge23}. A numerical analysis detects only one root $\bar{\lambda}$ of $g_{2,3}$ in the interval $0<\lambda<1$, and shows that the factors of $N_{2,3}$ do not have roots in this interval. The root $\bar{\lambda}$ is approximately given by $\bar{\lambda}\approx0.645232$, which yields $k(\bar{\lambda})\approx (0.213451,\,0.239928,\,0.947034)$.
    \end{itemize}
    
This proves items (i) (b), (i) (c), (i) (d), except for the instability property.

Finally, we continue with the interior $\mathrm{int}(\mathcal{R})$. In view of \eqref{reg red ico ham in the inv}, the critical points in this region are the real solutions $(s,t)\in U$ of the system:
\begin{equation*}
\begin{split}
\frac{\partial}{\partial s}p^I_{\mathbb{I}}(s,t)=
c F_4(s,t)^3 F_5(s,t)\,G_2(s,t)=0\\
\frac{\partial}{\partial t}p^I_{\mathbb{I}}(s,t)=
c F_4(s,t)^3 F_5(s,t)\,H_2(s,t)=0
\end{split}
\end{equation*} where the complete expressions of the polynomials $G_2$, $H_2$ are given in \eqref{icodsfactor} and \eqref{icodtfactor}. The $0$ level set of $p_{\mathbb{I}}$ coincides with the $0$ level set of $\tilde{h}_{\mathbb{I}}$ which contains only collision equilibria. There is no collision equilibrium in $\mathrm{int}(\mathcal{R})$, therefore $F_4$, $F_5$ do not simultaneously vanish in $U$ and can be ignored. We are left with two polynomials $G_2(s,t)$, $H_2(s,t)$ which might contain common real zeroes in $U$. Therefore, the problem of finding critical points of $p_{\mathbb{I}}$ in $\mathrm{int}(\mathcal{R})$ is reduced to solving $G_2(s,t)=0$, $H_2(s,t)=0$ for $(s,t)\in U\subset (2/15,\,4/27)\times (1/15,\,25/243)$ (see Remark \ref{rem: ico inv fund reg bounds}).
The $s$- and $t$-values of the relevant common real roots of $G_2$, $H_2$ are the real roots in these intervals of the resultants $\mathrm{Res}_t(G_2, H_2)$ and $\mathrm{Res}_s(G_2, H_2)$, respectively. These resultants allow a factorisation of the form:
\[\mathrm{Res}_s(G_2,H_2)(t)=q_4(t)\cdot L_4(t)\,,\]
\[\mathrm{Res}_t(G_2,H_2)(s)=q_3(s)\cdot L_3(s)\,.\]
The complete expressions of these resultants and of the polynomials $q_3(s)$ and $q_4(t)$ are given in \eqref{icointeqpoly2}, \eqref{icoresultants} and  \eqref{icointeqpoly1}, \eqref{icoresultantt}. The expressions for the polynomials $L_3(s)
$ and $L_4(
t)$ can easily be read from \eqref{icoresultants} and \eqref{icoresultantt}.

A numerical investigation reveals the existence of two zeroes $\bar{t}<\bar{t}'$ of $q_4$ in the interval $(1/15,\,25/243)$ and of a zero $\bar{s}$ of $q_3$ in the interval $(2/15,\,4/27)$. It can be checked (again, numerically) that the roots of the factors $L_4(t)$ and $L_3(s)$ do not lead to any common zero $(s,t)\in U$ of $G_2$, $H_2$ and only $(\bar{s},\bar{t})$ is actually a common zero of $G_2$ and $H_2$ in this region (see Figure \ref{fig: FullIcoInvReg}). The point $(\bar{s},\bar{t})$ in the invariant plane is approximately $(0.143135,\,0.091740)$. 

\begin{figure}
   \centering
   \begin{overpic}[width=0.98\textwidth]{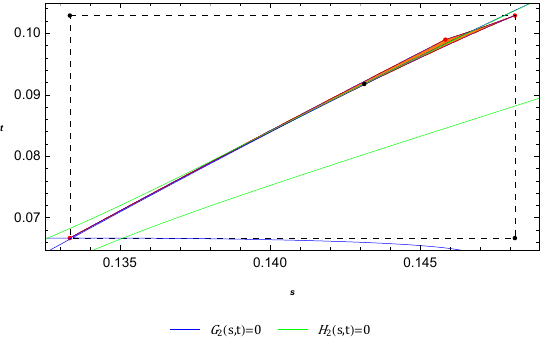}
    \put(10,19){\textcolor{red}{$V_
    3$}}
    \put(80,59){\textcolor{red}{$V_1$}}
    \put(96,60){\textcolor{red}{$V_2$}}
    \put(81, 22){\textcolor{black}{$(4/27,\,1/15)$}}
    \put(14,59){\textcolor{black}{$(2/15,\,25/243)$}}
    \put(64,45){\textcolor{black}{$(\bar{s},\bar{t})$}}
    \end{overpic}
    \caption{In orange, the domain $U=\psi^{-1}(\mathrm{int}(\mathcal{R}))$ in the invariants $(s, t)$-plane that parametrises the interior $\mathrm{int}(\mathcal{R})$ of the fundamental region. In blue and green the local vanishing loci of $G_2$, $H_2$.}
    \label{fig: FullIcoInvReg}
\end{figure}
In order to determine the (Lyapunov) stability of the equilibria found, we parametrise the northern hemisphere (that contains the fundamental region $\mathcal{R}$) by setting $z=\sqrt{1-x^2-y^2}$ and evaluate the Hessian of $Q(x, y)\coloneqq p_{\mathbb{I}}(x, y, \sqrt{1-x^2-y^2})$ at precise enough numerical approximations\footnote{In fact, if the Hessian is positive/negative definite or indefinite at a point $p\in\R^2$, it will be such at any point $p'$ in an open neighbourhood of $p$.} of the critical points in $\mathcal{R}$ (excluding $V_1$, $V_2$ and $V_3$); the remaining critical points have the same stability/instability properties by $\mathbb{I}_h$-symmetry. We have:
\begin{equation*}
\begin{aligned} &\mathrm{Hess}(Q)(\xi_1,\xi_2)\approx\begin{pmatrix}-3.787570\cdot 10^{23}&5.126704\cdot 10^{23}\\5.126704\cdot 10^{23}&-1.165732\cdot 10^{24}\end{pmatrix}\quad\text{Negative definite}\,,\\&\mathrm{Hess}(Q)(\sqrt{\zeta},0)\approx\begin{pmatrix}-1.356064\cdot 10^{24}&0\\0&6.196020\cdot 10^{23}\end{pmatrix}\quad\text{Indefinite}\,,\\&\mathrm{Hess}(Q)(0,\sqrt{\eta})\approx\begin{pmatrix}2.724941\cdot 10^{23}&0\\0&-1.704172\cdot 10^{24}\end{pmatrix}\quad\text{Indefinite}\,,\\ &\mathrm{Hess}(Q)(k(\bar{\lambda})_{x,y})\approx\begin{pmatrix}-3.886003\cdot 10^{23}&7.433536\cdot 10^{23}\\7.433536\cdot 10^{23}&-1.127358\cdot 10^{24}\end{pmatrix}\quad\text{Indefinite}\,,\end{aligned}
\end{equation*} which ends the proof.
\end{proof}

To complete the classification of icosahedrally symmetric equilibrium configurations, we would have to prove that the equilibria of the regularised reduced system \eqref{regredsys} described in item (i) of Theorem \ref{thm: Ico eqs} are the unique non-collision equilibria. Despite the simplifications introduced by the invariant coordinates, the degree and algebraic complexity of the resultants and polynomials appearing in the proof of Theorem \ref{thm: Ico eqs} remain high. Direct estimates become quite complicated, so we limit ourselves to conjecturing the following:
\begin{namedthm*}{\textbf{Conjecture}} The equilibria described in item (i) of Theorem \ref{thm: Ico eqs} are the only non-collision equilibria of the regularised reduced system \eqref{regredsys} (with $N=60$, $K=\mathbb{I}$) in the fundamental region $\mathcal{R}$.
\end{namedthm*} There are two main reasons which support our conjecture:
\begin{itemize}
\item The numerical investigation\footnote{Using \textit{Mathematica}'s $\mathrm{NSolve}$ with high working precision.} did not produce any evidence of the existence of other common real roots besides $(\bar{s},\bar{t})$ of the polynomials $G_2, H_2$ in $U$; nor did it show the existence of roots of the derivatives of the polynomials $p_{\mathbb{I}}(x,\,0,\,\sqrt{1-x^2})$, $p_{\mathbb{I}}(0,\,y,\,\sqrt{1-y^2})$ and $p_{\mathbb{I}}(k(\lambda))$ in the respective intervals except for $\sqrt{\zeta}$, $\sqrt{\eta}$ and $\bar{\lambda}$, respectively.

\item The number and type of equilibria of the regularised reduced system \eqref{regredsys} described in Theorem \ref{thm: Ico eqs} are consistent with the Poincaré--Hopf Index Theorem. In fact, the index of stable non-collision equilibria and collision equilibria is $+1$, while the index of unstable non-collision equilibria is $-1$, and we have that: \begin{equation*}\begin{aligned} &\#\text{collision eqs.}+\#\text{stable non-collision eqs.}-\#\text{unstable non-collision eqs.}=\\&=62+120-180=2=\chi(\mathbb{S}^2)\end{aligned}\end{equation*} where we have denoted by $\chi(\mathbb{S}^2)$ the Euler--Poincaré characteristic of $\mathbb{S}^2$.
\end{itemize}

\begin{figure}\centering
\begin{overpic}[width=.5\textwidth]{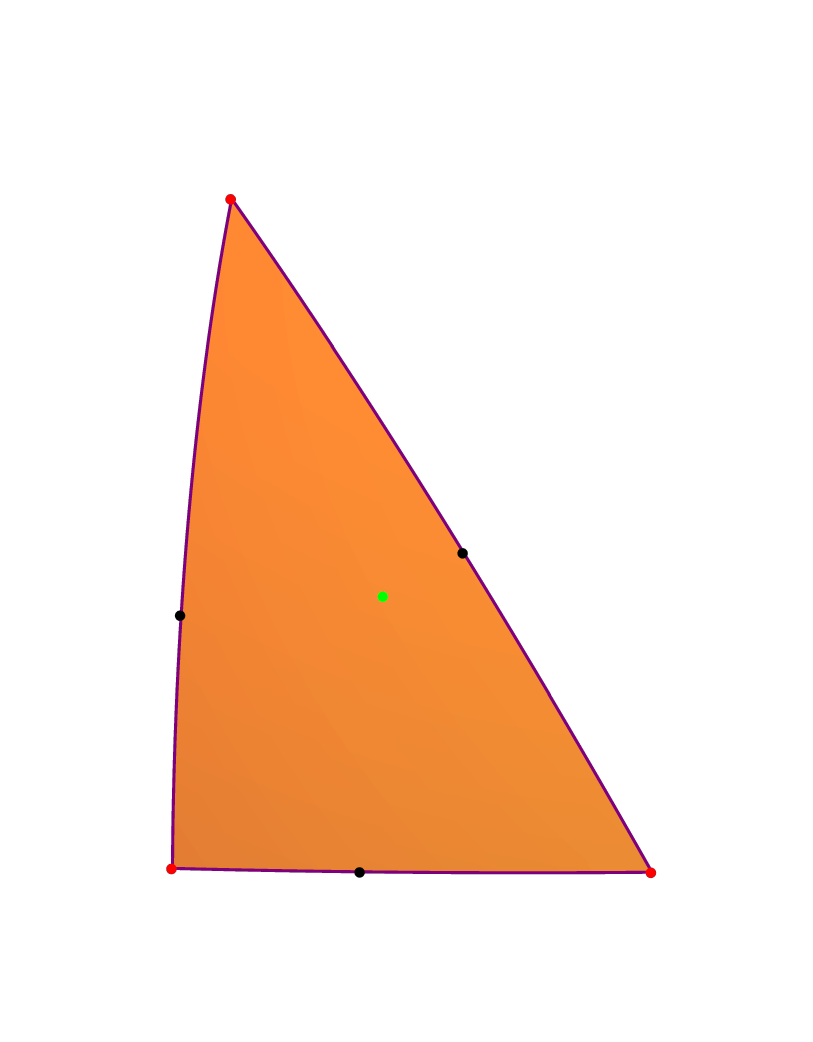}
    \put(12,13){\fontsize{7}{8.4}\selectfont\textcolor{red}{$V_1$}}
    \put(64,13){\fontsize{7}{8.4}\selectfont\textcolor{red}{$V_2$}}
    \put(21,83){\fontsize{7}{8.4}\selectfont\textcolor{red}{$V_3$}}
    
    \put(-3,40){\parbox{20mm}{\fontsize{7}{8.4}\selectfont\centering Irreg. Trunc.\newline Icosahedron}}
    
    \put(46,46){\parbox{20mm}{\fontsize{7}{8.4}\selectfont\centering Pseudo\newline Rhombicosid.}}

    \put(26,35){\parbox{20mm}{\fontsize{7}{8.4}\selectfont\centering \textcolor{green}{Pseudo Snub}\newline \textcolor{green}{Dodecahedron (R)}}}
    
    \put(26,11){\parbox{20mm}{\fontsize{7}{8.4}\selectfont\centering Irreg. Trunc.\newline Dodecahedron}}

    \end{overpic}
  \caption{Equilibria of the regularised reduced system \eqref{regredsys} on the fundamental region $\mathcal{R}$ (cf. Theorem \ref{thm: Ico eqs}).}
    \label{fig: FullIcoregwitheqs}
\end{figure}

\subsection{Reduced phase portrait}\label{sec: ico dynamics, ph portrait}
\begin{figure}
    \centering
    \renewcommand{\figurename}{Figure}\includegraphics[width=.85\linewidth]{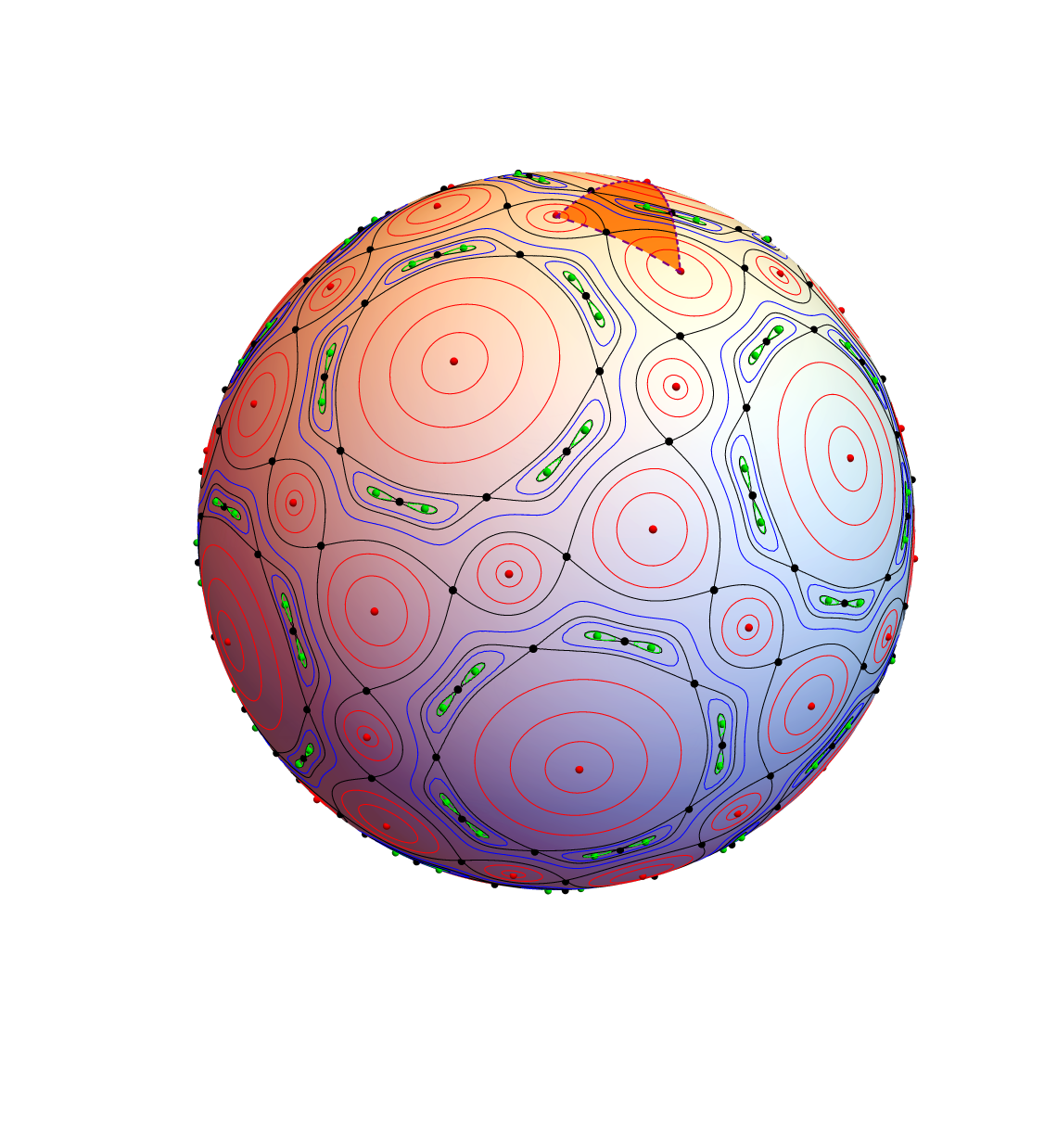}
    \caption{The (regularised) reduced phase portrait of icosahedrally symmetric solutions of the $60$-vortex problem on the sphere, with the fundamental region $\mathcal{R}$ highlighted.}
    \label{fig: IcoPhasePortrait}
\end{figure}

In analogy to the octahedral case, we present in Figure \ref{fig: IcoPhasePortrait} a numerical approximation of the phase portrait of the regularised reduced system \eqref{regredsys}. The pseudo snub dodecahedron equilibrium points are indicated in green, the collision equilibrium points (icosidodecahedron, dodecahedron and icosahedron) are indicated  in red and the truncated dodecahedron, truncated icosahedron and pseudo rhombicosidodecahedron equilibrium points are indicated in black. The colour code is for periodic orbits near the stable equilibria or heteroclinic/homoclinic orbits emanating from the unstable equilibria. Blue periodic orbits are other families of periodic orbits not emanating from equilibria. 

\section*{Acknowledgements}
\addcontentsline{toc}{section}{Acknowledgements} \markboth{Acknowledgements}{}
This article is based on the author's master's thesis \cite{Poma}. I am grateful to L. C. García-Naranjo for the many suggestions that helped improve the exposition.
\appendix
\section{Dimensions of polyhedra corresponding to equilibrium configurations}

In this appendix, we report the edge lengths and truncation ratios of the polyhedra whose vertices are occupied by the vortices in non-collision equilibrium configurations (of the full system \eqref{system}).

\subsection{Octahedral symmetry}
In the case of octahedrally symmetric non-collision equilibrium configurations (cf. item (i) of Theorem \ref{thm: Oct eqs}):
\begin{itemize}
    \item[(a)] The pseudo snub cube contains $6$ square faces with edges of length $\sqrt{2}\sqrt{\alpha_1^2+\alpha_2^2}\approx0.715734$ and $8$ equilateral triangular faces with edges of length \[\sqrt{(\alpha_1-\alpha_2)^2+(\alpha_1-\alpha_3)^2+(\alpha_2-\alpha_3)^2}\approx0.767210\] not adjacent to the square faces. The remaining  $24$ triangular faces are not regular, and have an edge in common with a square face and an edge in common with an equilateral triangular face; the third edge has length $\sqrt{4 \alpha_2^2 + 2 (\alpha_1 - \alpha_3)^2}\approx0.780569$.  
    
    \item[(b)] The truncated octahedron contains $6$ square faces with edge length $\sqrt{2}\sqrt{\beta}\approx0.671364$ and $8$ faces with the shape of an irregular hexagon (with edges of alternating length $\sqrt{2}\sqrt{\beta}\approx0.671364$ and $\sqrt{2} (\sqrt{1 - \beta} - \sqrt{\beta})\approx0.573333$). The truncation ratio is given by: \[R=\sqrt{\beta}/\big( \sqrt{1 - \beta} + \sqrt{\beta}\,\big)\approx0.350388\,.\]
    
    \item[(c)] The pseudo rhombicuboctahedron contains $6$ square faces with edge length $2\sqrt{\gamma}\approx0.700587$ and $8$ equilateral triangular faces with edge length $\sqrt{2} (\sqrt{1 - 2 \gamma} - \sqrt{\gamma})\approx0.733096$, not adjacent to the square faces. The remaining faces are rectangles (with edges of alternating length $2\sqrt{\gamma}\approx0.700587$ and $\sqrt{2} (\sqrt{1 - 2 \gamma} - \sqrt{\gamma})\approx0.733096$). 

    \item[(d)] The truncated cube contains $8$ equilateral triangular faces with edge length $\sqrt{2} (\sqrt{\delta}-\sqrt{1 - 2 \delta} )\approx0.623258$ and $6$ faces with the shape of  an irregular octagon (with edges of alternating length $\sqrt{2} (\sqrt{\delta}-\sqrt{1 - 2 \delta} )\approx0.623258$ and $2 \sqrt{1 - 2 \delta}\approx0.489739$). The truncation ratio is given by: \[R=\frac{1}{2}-\frac{1}{2}\sqrt{\frac{1-2\delta}{\delta}}\approx0.321414\,.\]
\end{itemize}

\subsection{Icosahedral symmetry}
In the case of icosahedrally symmetric non-collision equilibrium configurations (cf. item (i) of Theorem \ref{thm: Ico eqs}):

\begin{itemize}
    \item[(a)] The pseudo snub dodecahedron contains $12$ regular pentagonal faces with edges of length \[\sqrt{(3-\phi)\xi_1^2+(\xi_2+(1-\phi)\xi_3)^2}\approx0.431065\] and $20$ equilateral triangular faces with edges of length \[\sqrt{(\phi\,\xi_1+(1-\phi)\xi_3)^2+3\xi_2^2}\approx0.497250\] not adjacent to the pentagonal faces. The remaining  $60$ triangular faces are not regular, and have an edge in common with a pentagonal face and an edge in common with an equilateral triangular face; the third edge has length $2\sqrt{\xi_1^2+\xi_2^2}\approx0.517531$.  
    
    \item[(b)] The truncated dodecahedron contains $20$ equilateral triangular faces with edge length $(\phi-1)\sqrt{1-\zeta}-\phi\sqrt{\zeta}\approx0.377879$ and $12$ faces with the shape of an irregular decagon (with edges of alternating length $(\phi-1)\sqrt{1-\zeta}-\phi\sqrt{\zeta}\approx0.377879$ and $2\sqrt{\zeta}\approx0.288840$). The truncation ratio is given by: \[R=\frac{1}{2}-\frac{\phi}{2(\phi-1)}\sqrt{\frac{\zeta}{1-\zeta}}\approx0.308949\,.\]

    \item[(c)] The truncated icosahedron contains $12$ regular pentagonal faces with edge length $(\phi-1)\sqrt{1-\eta}-\sqrt{\eta}\approx0.415156$ and $20$ faces with the shape of an irregular hexagon (with edges of alternating length $(\phi-1)\sqrt{1-\eta}-\sqrt{\eta}\approx0.415156$ and $2\sqrt{\eta}\approx0.382892$). The truncation ratio is given by: \[R=\frac{1}{2}-\frac{1}{2(\phi-1)}\sqrt{\frac{\eta}{1-\eta}}\approx0.342198\,.\]
    
    \item[(d)] For brevity, we set $k(\bar{\lambda})=(\iota_1,\iota_2,\iota_3)$. The pseudo rhombicosidodecahedron contains $12$ regular pentagonal faces with edge length $2\iota_1\approx0.426902$ and $20$ equilateral triangular faces with edge length $2\iota_2\approx0.479856$, not adjacent to the pentagonal faces. The remaining $30$ faces are rectangles (with edges of alternating length $2\iota_1\approx0.426902$ and $2\iota_2\approx0.479856$).

\end{itemize}

\section{Complete expressions}

This appendix contains the complete expressions of Hamiltonians and other polynomials and functions appearing in previous sections.

\subsection{Octahedral case}

{\fontsize{7pt}{8.4pt}\selectfont
\begin{flalign}\label{octregham}
\begin{aligned}&\tilde{h}_{\mathbb{O}}(x,y,z)=\\
&\Big[2^{26}(x^2 + y^2)^3  (x^2 + z^2)^3  (y^2 + z^2)^3  (x^2 + 2 y^2 - 2 x z + 
   z^2)  (x^2 + 2 y^2 + 2 x z + z^2)  \\ 
    & \cdot(2 x^2 + y^2 - 2 y z + 
   z^2) (x^2 - x y + y^2 - x z - y z + z^2)^2  (x^2 + x y + y^2 + 
   x z - y z + z^2)^2 \\ 
   &\cdot(x^2 + x y + y^2 - x z + y z + z^2)^2 (x^2 - 
   x y + y^2 + x z + y z + z^2)^2  (2 x^2 + y^2 + 2 y z + 
   z^2)\\  
   &\cdot(x^2 - 2 x y + y^2 + 2 z^2)  (x^2 + 2 x y + y^2 + 2 z^2)\Big]^{12} 
\end{aligned}&&
\end{flalign}

\begin{flalign}\label{octsimplham}
\begin{aligned} &p_{\mathbb{O}}(x,y,z)=\\
&(x^2 + y^2)^3  (x^2 + z^2)^3  (y^2 + z^2)^3  (x^2 + 2 y^2 - 2 x z + 
   z^2)  (x^2 + 2 y^2 + 2 x z + z^2)  \\ 
     &\cdot(2 x^2 + y^2 - 2 y z + 
   z^2) (x^2 - x y + y^2 - x z - y z + z^2)^2  (x^2 + x y + y^2 + 
   x z - y z + z^2)^2 \\ 
   &\cdot(x^2 + x y + y^2 - x z + y z + z^2)^2 (x^2 - 
   x y + y^2 + x z + y z + z^2)^2  (2 x^2 + y^2 + 2 y z + 
   z^2)\\  
   &\cdot(x^2 - 2 x y + y^2 + 2 z^2)  (x^2 + 2 x y + y^2 + 2 z^2)
\end{aligned}&&
\end{flalign}

\begin{flalign}\label{octedgeqpoly1}\begin{aligned} &\frac{d}{dx}p_{\mathbb{O}}(x,0,\sqrt{1-x^2})=\\&2 (-1 + x)^2 x^5 (1 + x)^2 (-2 + x^2) (1 + x^2) (-1 + 2 x^2)(-1 + 
   x \sqrt{1 - x^2})^3 (1 + x\sqrt{1 - x^2})^3\\&  \cdot(-6 + 
   41 x^2 - 72 x^4 + 22 x^6 + 89 x^8 - 120 x^{10} + 40 x^{12})\end{aligned}&&
\end{flalign}

\begin{flalign}\label{octedgeqpoly2}\begin{aligned} &\frac{d}{dx}p_{\mathbb{O}}(x,x,\sqrt{1-2x^2})=\frac{d}{dx}p_{\mathbb{O}}(x,\sqrt{1-2x^2},x)=\\&-64 (-1 + x)^5 x^5 (1 + x)^5 (1 + x^2)^3 (-1 + 3 x^2) (-1 + x^2 - 
   2 x \sqrt{1 - 2 x^2})\\&\cdot (1 + x^2 - 2 x \sqrt{1 - 2 x^2}) (-1 + x^2 + 
   2 x \sqrt{1 - 2 x^2}) (1 + x^2 + 2 x \sqrt{1 - 2 x^2})\\&\cdot(-3 + 41 x^2 - 
   190 x^4 + 502 x^6 - 323 x^8 - 1023 x^{10} + 1188 x^{12})\end{aligned}&&
\end{flalign}

\begin{flalign}\label{octdsfactor}\begin{aligned}&G_1(s,t)=\\& +12 - 88 s + 97 s^2 + 90 s^3 - 71 s^4 - 40 s^5 + 316 t - 1282 s t + 
 2940 s^2 t - 382 s^3 t+ 28 s^4 t\\ & - 4743 t^2 + 3114 s t^2 + 
 117 s^2 t^2 - 216 t^3 + 108 s t^3
\end{aligned}&&
\end{flalign}

\begin{flalign}\label{octdtfactor}\begin{aligned}&H_1(s,t)=\\&-12 + 108 s - 129 s^2 + 642 s^3 + 27 s^4 + 12 s^5 - 336 t - 4446 s t+ 1152 s^2 t + 66 s^3 t + 10503 t^2\\ & - 5454 s t^2 + 
 135 s^2 t^2 - 2268 t^3
    \end{aligned}&&
\end{flalign}

\begin{flalign}\label{octinteqpoly1}\begin{aligned}&q_1(s)=\\&+7355241 - 49011072 s + 71977588 s^2 - 52132560 s^3 + 34212078 s^4  
 -5786528 s^5 + 14244 s^6\\ &  + 26736 s^7 + 529 s^8
\end{aligned}&&
\end{flalign}

\begin{flalign}\label{octinteqpoly2}\begin{aligned} &q_2(t)=\\&+44136711 - 4537885428  t - 37526316256  t^2 - 54799208564  t^3 + 
 316163424918  t^4\\& + 15809108548  t^5 - 7863707368  t^6 - 
 400482396  t^7 + 2299563  t^8\end{aligned} &&
\end{flalign}

\begin{flalign}\label{octresultants}\begin{aligned} &\mathrm{Res}_s(G_1,H_1)(t)=\\&78364164096 (-1 + t) (9 + t) (9 - 294 t + t^2)(1 - 14 t + t^2) \\&\cdot(1 + 
   34 t + t^2) (3 + 14 t + 27  t^2)\,q_2(t)
\end{aligned}&&\end{flalign}

\begin{flalign}\label{octresultantt}\begin{aligned} &\mathrm{Res}_t(G_1,H_1)(s)=\\&5038848 (-7 + s) (-1 + s) (1 - 14 s + s^2)\\ &\cdot (1 + 34 s + s^2) (-23 + 
   58 s + s^2) (11 - 2 s + 3  s^2)\, q_1(s)
    \end{aligned} &&
\end{flalign}
}

\subsection{Icosahedral case}

{\fontsize{7pt}{8.4pt}\selectfont
\begin{flalign}\label{icoregham}\begin{aligned} &\tilde{h}_{\mathbb{I}}(x,y,z)=\\&\Big[\frac{1}{2^{34}}(x^2 + y^2) (-3 x^2 + \sqrt{5} x^2 - 4 x y - 3 y^2 - \sqrt{5} y^2 - 
   6 z^2)^2 (-3 x^2 + \sqrt{5} x^2 + 4 x y - 3 y^2 - \sqrt{5} y^2 - 
   6 z^2)^2\\& \cdot (-3 x^2 + \sqrt{5} x^2 - 5 y^2 + \sqrt{5} y^2 + 2 x z - 
   2 \sqrt{5} x z - 2 z^2)^2 (-3 x^2 + \sqrt{5} x^2 - 5 y^2 + 
   \sqrt{5} y^2 - 2 x z + 2 \sqrt{5} x z - 2 z^2)^2 \\& \cdot(x^2 + z^2) (y^2 + 
   z^2) (x^2 - x y + y^2 - x z - y z + z^2)^2 (x^2 + x y + y^2 + x z -
    y z + z^2)^2 \\&  \cdot(x^2 + x y + y^2 - x z + y z + z^2)^2(x^2 - x y + 
   y^2 + x z + y z + z^2)^2 \\&  \cdot(5 x^2 + \sqrt{5} x^2 + 3 y^2 + 
   \sqrt{5} y^2 - 2 y z - 2 \sqrt{5} y z + 2 z^2)^2(5 x^2 + 
   \sqrt{5} x^2 + 3 y^2 + \sqrt{5} y^2 + 2 y z + 2 \sqrt{5} y z + 
   2 z^2)^2\\&  \cdot(5 x^2 + \sqrt{5} x^2 - 4 x y + 5 y^2 - \sqrt{5} y^2 + 
   2 x z - 2 \sqrt{5} x z - 2 y z - 2 \sqrt{5} y z + 6 z^2)\\& \cdot (5 x^2 + 
   \sqrt{5} x^2 + 4 x y + 5 y^2 - \sqrt{5} y^2 - 2 x z + 
   2 \sqrt{5} x z - 2 y z - 2 \sqrt{5} y z + 6 z^2)\\& \cdot (5 x^2 + 
   \sqrt{5} x^2 + 4 x y + 5 y^2 - \sqrt{5} y^2 + 2 x z - 
   2 \sqrt{5} x z + 2 y z + 2 \sqrt{5} y z + 6 z^2)\\& \cdot (5 x^2 + 
   \sqrt{5} x^2 - 4 x y + 5 y^2 - \sqrt{5} y^2 - 2 x z + 
   2 \sqrt{5} x z + 2 y z + 2 \sqrt{5} y z + 6 z^2)\\& \cdot (-5 x^2 + 
   \sqrt{5} x^2 + 2 x y + 2 \sqrt{5} x y - 6 y^2 - 4 x z + 2 y z - 
   2 \sqrt{5} y z - 5 z^2 - \sqrt{5} z^2)\\& \cdot (-5 x^2 + \sqrt{5} x^2 - 
   2 x y - 2 \sqrt{5} x y - 6 y^2 + 4 x z + 2 y z - 2 \sqrt{5} y z - 
   5 z^2 - \sqrt{5} z^2)\\& \cdot (-5 x^2 + \sqrt{5} x^2 - 2 x y - 
   2 \sqrt{5} x y - 6 y^2 - 4 x z - 2 y z + 2 \sqrt{5} y z - 5 z^2 - 
   \sqrt{5} z^2) \\& \cdot(-5 x^2 + \sqrt{5} x^2 + 2 x y + 2 \sqrt{5} x y - 
   6 y^2 + 4 x z - 2 y z + 2 \sqrt{5} y z - 5 z^2 - 
   \sqrt{5} z^2)\\& \cdot (3 x^2 + \sqrt{5} x^2 + 6 y^2 - 4 x z + 3 z^2 - 
   \sqrt{5} z^2)^2  (3 x^2 + \sqrt{5} x^2 + 6 y^2 + 4 x z + 3 z^2 - 
   \sqrt{5}z^2)^2 \\& \cdot(2 x^2 + 2 x y - 2 \sqrt{5} x y + 3 y^2 - 
   \sqrt{5} y^2 + 5 z^2 - \sqrt{5} z^2)^2 (2 x^2 - 2 x y + 
   2 \sqrt{5} x y + 3 y^2 - \sqrt{5} y^2 + 5 z^2 - 
  \sqrt{5} z^2)^2 \\& \cdot(6 x^2 + 2 x y - 2 \sqrt{5} x y + 5 y^2 + 
   \sqrt{5}y^2 - 2 x z - 2 \sqrt{5}x z - 4 y z + 5 z^2 - 
   \sqrt{5} z^2)\\& \cdot (6 x^2 - 2 x y + 2 \sqrt{5} x y + 5 y^2 + 
   \sqrt{5} y^2 + 2 x z + 2 \sqrt{5} x z - 4 y z + 5 z^2 - 
   \sqrt{5} z^2)\\& \cdot (6 x^2 - 2 x y + 2 \sqrt{5} x y + 5 y^2 + 
   \sqrt{5} y^2 - 2 x z - 2 \sqrt{5} x z + 4 y z + 5 z^2 - 
   \sqrt{5} z^2)\\& \cdot (6 x^2 + 2 x y - 2 \sqrt{5}x y + 5 y^2 + 
   \sqrt{5} y^2 + 2 x z + 2 \sqrt{5} x z + 4 y z + 5 z^2 - 
   \sqrt{5}z^2)\\& \cdot (-5 x^2 + \sqrt{5} x^2 - 2 y^2 + 2 y z - 
   2 \sqrt{5} y z - 3 z^2 + \sqrt{5} z^2)^2  (-5 x^2 + \sqrt{5} x^2 - 
   2 y^2 - 2 y z + 2 \sqrt{5} y z - 3 z^2 + \sqrt{5} z^2)^2\\& \cdot(2 x^2 + 
   5 y^2 + \sqrt{5} y^2 - 2 x z - 2 \sqrt{5} x z + 3 z^2 + 
   \sqrt{5} z^2)^2  (2 x^2 + 5 y^2 + \sqrt{5} y^2 + 2 x z + 
   2 \sqrt{5} x z + 3 z^2 + \sqrt{5} z^2)^2\\&  \cdot(6 x^2 + 3 y^2 - 
   \sqrt{5} y^2 - 4 y z + 3 z^2 + \sqrt{5} z^2)^2  (6 x^2 + 3 y^2 - 
   \sqrt{5} y^2 + 4 y z + 3 z^2 + \sqrt{5} z^2)^2\\& \cdot (3 x^2 + 
   \sqrt{5} x^2 - 2 x y - 2 \sqrt{5}x y + 2 y^2 + 5 z^2 + 
  \sqrt{5} z^2)^2 (3 x^2 + \sqrt{5} x^2 + 2 x y + 2 \sqrt{5} x y + 
   2 y^2 + 5 z^2 + \sqrt{5} z^2)^2\Big]^{30}\end{aligned}&&\end{flalign}

\begin{flalign}\label{icosimplregham}\begin{aligned} &p_{\mathbb{I}}(x,y,z)=\\&(x^2 + y^2) (-3 x^2 + \sqrt{5} x^2 - 4 x y - 3 y^2 - \sqrt{5} y^2 - 
   6 z^2)^2 (-3 x^2 + \sqrt{5} x^2 + 4 x y - 3 y^2 - \sqrt{5} y^2 - 
   6 z^2)^2\\& \cdot (-3 x^2 + \sqrt{5} x^2 - 5 y^2 + \sqrt{5} y^2 + 2 x z - 
   2 \sqrt{5} x z - 2 z^2)^2 (-3 x^2 + \sqrt{5} x^2 - 5 y^2 + 
   \sqrt{5} y^2 - 2 x z + 2 \sqrt{5} x z - 2 z^2)^2 \\& \cdot(x^2 + z^2) (y^2 + 
   z^2) (x^2 - x y + y^2 - x z - y z + z^2)^2 (x^2 + x y + y^2 + x z -
    y z + z^2)^2\\& \cdot (x^2 + x y + y^2 - x z + y z + z^2)^2(x^2 - x y + 
   y^2 + x z + y z + z^2)^2\\& \cdot (5 x^2 + \sqrt{5} x^2 + 3 y^2 + 
   \sqrt{5} y^2 - 2 y z - 2 \sqrt{5} y z + 2 z^2)^2(5 x^2 + 
   \sqrt{5} x^2 + 3 y^2 + \sqrt{5} y^2 + 2 y z + 2 \sqrt{5} y z + 
   2 z^2)^2\\&  \cdot(5 x^2 + \sqrt{5} x^2 - 4 x y + 5 y^2 - \sqrt{5} y^2 + 
   2 x z - 2 \sqrt{5} x z - 2 y z - 2 \sqrt{5} y z + 6 z^2)\\& \cdot (5 x^2 + 
   \sqrt{5} x^2 + 4 x y + 5 y^2 - \sqrt{5} y^2 - 2 x z + 
   2 \sqrt{5} x z - 2 y z - 2 \sqrt{5} y z + 6 z^2)\\& \cdot (5 x^2 + 
   \sqrt{5} x^2 + 4 x y + 5 y^2 - \sqrt{5} y^2 + 2 x z - 
   2 \sqrt{5} x z + 2 y z + 2 \sqrt{5} y z + 6 z^2)\\& \cdot (5 x^2 + 
   \sqrt{5} x^2 - 4 x y + 5 y^2 - \sqrt{5} y^2 - 2 x z + 
   2 \sqrt{5} x z + 2 y z + 2 \sqrt{5} y z + 6 z^2)\\& \cdot (-5 x^2 + 
   \sqrt{5} x^2 + 2 x y + 2 \sqrt{5} x y - 6 y^2 - 4 x z + 2 y z - 
   2 \sqrt{5} y z - 5 z^2 - \sqrt{5} z^2)\\& \cdot (-5 x^2 + \sqrt{5} x^2 - 
   2 x y - 2 \sqrt{5} x y - 6 y^2 + 4 x z + 2 y z - 2 \sqrt{5} y z - 
   5 z^2 - \sqrt{5} z^2)\\& \cdot (-5 x^2 + \sqrt{5} x^2 - 2 x y - 
   2 \sqrt{5} x y - 6 y^2 - 4 x z - 2 y z + 2 \sqrt{5} y z - 5 z^2 - 
   \sqrt{5} z^2) \\& \cdot(-5 x^2 + \sqrt{5} x^2 + 2 x y + 2 \sqrt{5} x y - 
   6 y^2 + 4 x z - 2 y z + 2 \sqrt{5} y z - 5 z^2 - 
   \sqrt{5} z^2)\\& \cdot (3 x^2 + \sqrt{5} x^2 + 6 y^2 - 4 x z + 3 z^2 - 
   \sqrt{5} z^2)^2  (3 x^2 + \sqrt{5} x^2 + 6 y^2 + 4 x z + 3 z^2 - 
   \sqrt{5}z^2)^2 \\& \cdot(2 x^2 + 2 x y - 2 \sqrt{5} x y + 3 y^2 - 
   \sqrt{5} y^2 + 5 z^2 - \sqrt{5} z^2)^2 (2 x^2 - 2 x y + 
   2 \sqrt{5} x y + 3 y^2 - \sqrt{5} y^2 + 5 z^2 - 
  \sqrt{5} z^2)^2 \\& \cdot(6 x^2 + 2 x y - 2 \sqrt{5} x y + 5 y^2 + 
   \sqrt{5}y^2 - 2 x z - 2 \sqrt{5}x z - 4 y z + 5 z^2 - 
   \sqrt{5} z^2)\\& \cdot (6 x^2 - 2 x y + 2 \sqrt{5} x y + 5 y^2 + 
   \sqrt{5} y^2 + 2 x z + 2 \sqrt{5} x z - 4 y z + 5 z^2 - 
   \sqrt{5} z^2)\\& \cdot (6 x^2 - 2 x y + 2 \sqrt{5} x y + 5 y^2 + 
   \sqrt{5} y^2 - 2 x z - 2 \sqrt{5} x z + 4 y z + 5 z^2 - 
   \sqrt{5} z^2)\\& \cdot (6 x^2 + 2 x y - 2 \sqrt{5}x y + 5 y^2 + 
   \sqrt{5} y^2 + 2 x z + 2 \sqrt{5} x z + 4 y z + 5 z^2 - 
   \sqrt{5}z^2)\\& \cdot (-5 x^2 + \sqrt{5} x^2 - 2 y^2 + 2 y z - 
   2 \sqrt{5} y z - 3 z^2 + \sqrt{5} z^2)^2  (-5 x^2 + \sqrt{5} x^2 - 
   2 y^2 - 2 y z + 2 \sqrt{5} y z - 3 z^2 + \sqrt{5} z^2)^2\\& \cdot(2 x^2 + 
   5 y^2 + \sqrt{5} y^2 - 2 x z - 2 \sqrt{5} x z + 3 z^2 + 
   \sqrt{5} z^2)^2  (2 x^2 + 5 y^2 + \sqrt{5} y^2 + 2 x z + 
   2 \sqrt{5} x z + 3 z^2 + \sqrt{5} z^2)^2\\&  \cdot(6 x^2 + 3 y^2 - 
   \sqrt{5} y^2 - 4 y z + 3 z^2 + \sqrt{5} z^2)^2  (6 x^2 + 3 y^2 - 
   \sqrt{5} y^2 + 4 y z + 3 z^2 + \sqrt{5} z^2)^2\\& \cdot (3 x^2 + 
   \sqrt{5} x^2 - 2 x y - 2 \sqrt{5}x y + 2 y^2 + 5 z^2 + 
  \sqrt{5} z^2)^2 (3 x^2 + \sqrt{5} x^2 + 2 x y + 2 \sqrt{5} x y + 
   2 y^2 + 5 z^2 + \sqrt{5} z^2)^2\end{aligned}&&\end{flalign}

\begin{flalign}\label{icoinv6}\begin{aligned} &\theta_2(x,y,z)=\\&\frac{1}{96} \Big(14 x^6 + 14 y^6 + 3 (13 + \sqrt{5}) y^4 z^2 - 
   3 (-13 + \sqrt{5}) y^2 z^4 + 14 z^6 + 
   3 x^4 ((13 + \sqrt{5}) y^2\\& - (-13 + \sqrt{5}) z^2) + 
   3 x^2 (-((-13 + \sqrt{5}) y^4) + 36 y^2 z^2 + (13 + \sqrt{5}) z^4)\Big)
\end{aligned}&&\end{flalign}

\begin{flalign}\label{icoinv10}\begin{aligned} &\theta_3(x,y,z)=\\&\frac{1}{384} \Big(38 x^{10} + 38 y^{10} + 9 (17 + 3 \sqrt{5}) y^8 z^2 + 
   42 (7 + \sqrt{5}) y^6 z^4 - 42 (-7 + \sqrt{5}) y^4 z^6 + 
   9 (17 - 3 \sqrt{5}) y^2 z^8\\& + 38 z^{10} + 
   9 x^8 ((17 + 3 \sqrt{5}) y^2 + (17 - 3 \sqrt{5}) z^2) + 
   42 x^6 ((7 + \sqrt{5}) y^4 + 24 y^2 z^2 - (-7 + \sqrt{5}) z^4)\\& - 
   42 x^4 ((-7 + \sqrt{5}) y^6- 30 y^4 z^2 - 
      30 y^2 z^4 - (7 + \sqrt{5}) z^6) + 
   9 x^2 ((17 - 3 \sqrt{5}) y^8 + 112 y^6 z^2 + 140 y^4 z^4\\& + 
      112 y^2 z^6 + (17 + 3 \sqrt{5}) z^8)\Big)
\end{aligned}&&\end{flalign}

\begin{flalign}\label{Icopolyedge12}\begin{aligned} &g_{1,2}(u)= \\&-6750 + (249750 + 53100 \sqrt{5}) u + (-145600 - 
    911550 \sqrt{5}) u^2 + (-6421890 + 
    3252906 \sqrt{5}) u^3 \\& + (44430195 - 
    2092435 \sqrt{5}) u^4 + (-171529875 - 
    12520449 \sqrt{5}) u^5 + (418546015 + 
    35998909 \sqrt{5}) u^6 \\&+ (-675586370 - 
    45620334 \sqrt{5}) u^7 + (737129985 + 
    33819213 \sqrt{5}) u^8 + (-535576615 - 
    22301115 \sqrt{5}) u^9 \\&+ (240896645 + 
    26498965 \sqrt{5}) u^{10} + (-52378480 - 
    30705252 \sqrt{5}) u^{11} + (-4617080 + 
    20687642 \sqrt{5}) u^{12}\\& + (7561570 - 
    7186200 \sqrt{5}) u^{13} + (-2936250 + 1026600 \sqrt{5}) u^{14} + 
 391500 u^{15}
\end{aligned}&&\end{flalign}

\begin{flalign}\label{Icopolyedge13}\begin{aligned} &g_{1,3}(v)=\\&+6750 + (-249750 + 53100  \sqrt{5})  v + (145600 - 
     911550  \sqrt{5})  v^2 + (6421890 + 
     3252906  \sqrt{5})  v^3 \\&+ (-44430195 - 
     2092435  \sqrt{5})  v^4 + (171529875 - 
     12520449  \sqrt{5})  v^5 + (-418546015 + 
     35998909  \sqrt{5})  v^6\\& + (675586370 - 
     45620334  \sqrt{5})  v^7 + (-737129985 + 
     33819213  \sqrt{5})  v^8 + (535576615 - 
     22301115  \sqrt{5})  v^9\\& + (-240896645 + 
     26498965  \sqrt{5})  v^{10} + (52378480 - 
     30705252  \sqrt{5})  v^{11} + (4617080 + 
     20687642  \sqrt{5})  v^{12} \\&+ (-7561570 - 
     7186200 \sqrt{5})  v^{13} + (2936250 + 
     1026600  \sqrt{5})  v^{14} - 391500  v^{15}
\end{aligned}&&\end{flalign}

\begin{flalign}\label{Icopolyedge23}\begin{aligned}
&g_{2,3}(\lambda)=\\&+2396275200 + 
  1071646848 \sqrt{5}+ (-47925504000 - 21432936960 \sqrt{5}) \lambda + (462153344960 + 
     206681259072 \sqrt{5}) \lambda^2\\& + (-2860133672800 - 
     1279090663424 \sqrt{5}) \lambda^3 + (12755405292480 + 
     5704390662912 \sqrt{5}) \lambda^4\\& + (-43633637284840 - 
     19513555814904 \sqrt{5}) \lambda^5 + (118978356385920 + 
     53208738546028 \sqrt{5}) \lambda^6\\& + (-265277064270080 - 
     118635509715648\sqrt{5}) \lambda^7 + (492176881150660 + 
     220108192639004 \sqrt{5}) \lambda^8\\& + (-769092272400870 - 
     343948520397856\sqrt{5}) \lambda^9 + (1020316532581720 + 
     456299425021636 \sqrt{5}) \lambda^{10}\\& + (-1154335876741630 - 
     516234697620538 \sqrt{5}) \lambda^{11} + (1114792177317900 + 
     498550217121200 \sqrt{5}) \lambda^{12}\\& + (-916032463255115 - 
     409662169483383 \sqrt{5}) \lambda^{13} + (634350506331595 + 
     283690165948953 \sqrt{5}) \lambda^{14}\\& + (-362343822870780 - 
     162045073735168 \sqrt{5}) \lambda^{15} + (162246333721955 + 
     72558747485139 \sqrt{5}) \lambda^{16}\\& + (-48456351705540 - 
     21670308415588 \sqrt{5}) \lambda^{17} + (896440674460 + 
     400855663330 \sqrt{5}) \lambda^{18} \\&+ (10414394051730 + 
     4657515905898 \sqrt{5}) \lambda^{19} + (-8322015955110 - 
     3721782889970 \sqrt{5}) \lambda^{20}\\& + (4262750281790 + 
     1906422261744 \sqrt{5}) \lambda^{21} + (-1661732322955 - 
     743200827247 \sqrt{5}) \lambda^{22}\\& + (515463383310 + 
     230557224058 \sqrt{5}) \lambda^{23} + (-128536561860 - 
     57501566286 \sqrt{5}) \lambda^{24} \\&+ (25604101575 + 
     11456670875 \sqrt{5}) \lambda^{25} + (-3997471975 - 
     1788046725 \sqrt{5}) \lambda^{26}\\& + (473382000 + 
     210621600 \sqrt{5}) \lambda^{27} + (-40513500 - 
     17405100 \sqrt{5}) \lambda^{28} \\&+ (2430000 + 810000 \sqrt{5}) \lambda^{29} - 
  108000 \lambda^{30}\end{aligned}&&\end{flalign}

\begin{flalign}\label{Icopolyderivedge12}\begin{aligned} &\frac{d}{dx}p_{\mathbb{I}}(x,\,0,\,\sqrt{1-x^2})=\\&-\frac{4096}{25} x (5 + \sqrt{5} - 10 x^2)^3 (-3 + \sqrt{5} - 2 x^2)^4 (5 + 
   \sqrt{5} - 2 x^2)^4 (-3 + \sqrt{5} + 6 x^2) (-3 - \sqrt{5} - 3 x^2 + 
   \sqrt{5} x^2)^3\\& \cdot(5 - \sqrt{5} - 3 x^2 + \sqrt{5} x^2)^3 (-6 + 3 x^2 + 
   \sqrt{5} x^2)^3 (2 + 3 x^2 + \sqrt{5} x^2)^3 (-5 - \sqrt{5} + 
   2 \sqrt{5} x^2 - 4 x \sqrt{1 - x^2})\\& \cdot(3 - \sqrt{5} + 2 \sqrt{5} x^2 - 
   4 x \sqrt{1 - x^2}) (-1 + x \sqrt{1 - x^2})^3 (1 + 
   x \sqrt{1 - x^2})^3 (-5 - \sqrt{5} + 2 \sqrt{5} x^2 + 
   4 x \sqrt{1 - x^2})\\& \cdot(3 - \sqrt{5} + 2 \sqrt{5} x^2 + 
   4 x \sqrt{1 - x^2})(-3 - \sqrt{5} + x^2 + \sqrt{5} x^2 - 
   2 x \sqrt{1 - x^2} - 2 \sqrt{5} x \sqrt{1 - x^2})\\& \cdot(5 - \sqrt{5} + x^2 +
    \sqrt{5} x^2 - 2 x \sqrt{1 - x^2} - 2 \sqrt{5} x \sqrt{1 - x^2})(-2 -
    x^2 + \sqrt{5} x^2 + 2 x \sqrt{1 - x^2} - 
   2 \sqrt{5} x \sqrt{1 - x^2})\\& \cdot(6 - x^2 + \sqrt{5} x^2 + 
   2 x \sqrt{1 - x^2}- 2 \sqrt{5} x \sqrt{1 - x^2}) (-2 - x^2 + 
   \sqrt{5} x^2 - 2 x \sqrt{1 - x^2}+ 2 \sqrt{5} x \sqrt{1 - x^2})\\& \cdot(6 - 
   x^2 + \sqrt{5} x^2 - 2 x \sqrt{1 - x^2}+ 
   2 \sqrt{5} x \sqrt{1 - x^2}) (-3 - \sqrt{5} + x^2 + \sqrt{5} x^2 + 
   2 x \sqrt{1 - x^2} + 2 \sqrt{5} x \sqrt{1 - x^2}) \\& \cdot(5 - \sqrt{5} + x^2 +
    \sqrt{5} x^2 + 2 x \sqrt{1 - x^2}+ 2 \sqrt{5} x \sqrt{1 - x^2}) g_{1,2}(x^2)\end{aligned}&&\end{flalign}

\begin{flalign}\label{Icopolyderivedge13}\begin{aligned} &\frac{d}{dy}p_{\mathbb{I}}(0,\,y,\,\sqrt{1-y^2})=\\&-\frac{4096}{25} y (-5 + \sqrt{5} + 2 y^2)^4 (3 + \sqrt{5} + 2 y^2)^4 (-3 - 
\sqrt{5} + 6 y^2) (-5 + \sqrt{5} + 10 y^2)^3 (-2 - 3 y^2 + 
\sqrt{5} y^2)^3\\&\cdot (6 - 3 y^2 + \sqrt{5} y^2)^3 (-5 - \sqrt{5} + 3 y^2 + \sqrt{5} y^2)^3 (3 - \sqrt{5} + 3 y^2 + \sqrt{5} y^2)^3 (-3 - \sqrt{5} + 2 \sqrt{5} y^2 - 4 y \sqrt{1 - y^2})\\& \cdot(5 - \sqrt{5} + 
2 \sqrt{5} y^2 - 4 y \sqrt{1 - y^2}) (-1 + y \sqrt{1 - y^2})^3 (1 + y \sqrt{1 - y^2})^3 (-3 - \sqrt{5} + 2 \sqrt{5} y^2 + 
4 y \sqrt{1 - y^2}) \\& \cdot(5 - \sqrt{5} + 2 \sqrt{5} y^2 + 
4 y \sqrt{1 - y^2}) (-6 + y^2 + \sqrt{5} y^2 - 2 y \sqrt{1 - y^2} - 
2 \sqrt{5} y \sqrt{1 - y^2}) \\& \cdot(2 + y^2 + \sqrt{5} y^2 - 
2 y \sqrt{1 - y^2} - 2 \sqrt{5} y \sqrt{1 - y^2})(-5 - \sqrt{5} - 
y^2 + \sqrt{5} y^2 + 2 y \sqrt{1 - y^2} - 
2 \sqrt{5} y \sqrt{1 - y^2})\\& \cdot (3 - \sqrt{5} - y^2 + \sqrt{5} y^2 + 
2 y \sqrt{1 - y^2} - 2 \sqrt{5} y \sqrt{1 - y^2}) (-5 - \sqrt{5} - 
y^2 + \sqrt{5} y^2 - 2 y \sqrt{1 - y^2} + 
2 \sqrt{5} y \sqrt{1 - y^2})\\& \cdot(3 - \sqrt{5} - y^2 + \sqrt{5} y^2 - 
2 y \sqrt{1 - y^2} + 2 \sqrt{5} y \sqrt{1 - y^2}) (-6 + y^2 + 
\sqrt{5} y^2 + 2 y \sqrt{1 - y^2} + 2 \sqrt{5} y \sqrt{1 - y^2}) \\& \cdot(2 + y^2 + \sqrt{5} y^2 + 2 y \sqrt{1 - y^2} + 2 \sqrt{5} y \sqrt{1 - y^2}) g_{1,3}(y^2) \end{aligned}&&\end{flalign}

\begin{flalign}\label{Icopolyderivedge23}\begin{aligned} &\frac{d}{d\lambda}p_{\mathbb{I}}(k(\lambda))=\\&\frac{1}{(21 + 9 \sqrt{5} - 28 \lambda - 12 \sqrt{5} \lambda + 12 \lambda^2 + 4 \sqrt{5} \lambda^2)^{59}}\cdot \\&\Big[ -\left(\frac{1}{5}
  \right)(187072209578355573530071658587684226515959365500928) (3 + \sqrt{5})^3 (7 + 3 \sqrt{5})^2\\&\cdot (87403803 + 
   39088169 \sqrt{5}) (3 + \sqrt{5} - 4 \lambda) (1 + \sqrt{5} - 
   2 \lambda)^3 (3 + \sqrt{5} - 
   2 \lambda) (-1 + \lambda)^7 \lambda^3\\& \cdot (-10 - 4 \sqrt{5} + 
   15 \lambda + 5 \sqrt{5} \lambda - 10 \lambda^2) (-3 - \sqrt{5} + 3 \lambda + \sqrt{5} \lambda - 2 \lambda^2)^4\\& \cdot(-7 - 
   3 \sqrt{5} + 6 \lambda + 2 \sqrt{5} \lambda - 
   2 \lambda^2) (7 + 3 \sqrt{5} - 8 \lambda - 
   4 \sqrt{5} \lambda + 3 \lambda^2 + \sqrt{5} \lambda^2)\\& \cdot(7 + 
   3 \sqrt{5} - 7 \lambda - 3 \sqrt{5} \lambda + 3 \lambda^2 + 
   \sqrt{5} \lambda^2)^3 (14 + 6 \sqrt{5} - 14 \lambda - 
   6 \sqrt{5} \lambda + 5 \lambda^2 +\sqrt{5} \lambda^2) \\& \cdot(7 + 
   3 \sqrt{5} - 10 \lambda - 4 \sqrt{5} \lambda + 5 \lambda^2 + 
   \sqrt{5} \lambda^2)^3 (3 + \sqrt{5} - 6 \lambda - 
   2 \sqrt{5}\lambda + 5 \lambda^2 + \sqrt{5} \lambda^2) \\& \cdot(7 + 
   3 \sqrt{5} - 11 \lambda- 5 \sqrt{5}\lambda + 5 \lambda^2 + 
   2 \sqrt{5} \lambda^2) (14 + \sqrt{5} - 17 \lambda - 
   7 \sqrt{5} \lambda + 6 \lambda^2 + 
   2 \sqrt{5} \lambda^2)^3\\& \cdot(7 + 3 \sqrt{5} - 11 \lambda - 
   5 \sqrt{5} \lambda + 6 \lambda^2 + 
   2 \sqrt{5} \lambda^2)^3 (18 + 8 \sqrt{5} - 22 \lambda - 
   10 \sqrt{5} \lambda + 7 \lambda^2 + 
   3 \sqrt{5} \lambda^2)\\& \cdot(14 + 6 \sqrt{5} - 18 \lambda - 
   8 \sqrt{5} \lambda + 7 \lambda^2 + 
   3 \sqrt{5} \lambda^2)^3 (18 + 8 \sqrt{5} - 22 \lambda- 
   10 \sqrt{5} \lambda + 9 \lambda^2 + 
   3 \sqrt{5} \lambda^2)\\& \cdot(14 + 6 \sqrt{5} - 21 \lambda - 
   9 \sqrt{5} \lambda + 9 \lambda^2 + 
   3 \sqrt{5} \lambda^2)^3 (14 + 6 \sqrt{5} - 20 \lambda - 
   8 \sqrt{5} \lambda + 9 \lambda^2 + 3 \sqrt{5} \lambda^2)\\& \cdot(21 + 
   9 \sqrt{5} - 28 \lambda - 12 \sqrt{5}\lambda + 10 \lambda^2 + 
   4 \sqrt{5}\lambda^2) (18 + 8 \sqrt{5} - 25 \lambda - 
   11 \sqrt{5} \lambda + 10 \lambda^2 + 4 \sqrt{5} \lambda^2)^3 g_{2,3}(\lambda) \Big]\end{aligned}&&\end{flalign}
   
\begin{flalign}\label{icodsfactor}\begin{aligned} &G_2(s,t)=\\ & 6 \Big(-4174247 + 845119440 s - 20049967800 s^2 + 188979227280 s^3 - 
   861187991715 s^4 + 2013167663856 s^5\\ & - 2347664565150 s^6 + 
   209562932700 s^7 + 2994158277000 s^8 + 46747125000 s^9 - 
   767686095 t\\ & + 20135180340 s t - 181940184720 s^2 t + 
   692314639920 s^3 t - 1157398107525 s^4 t + 893590354080 s^5 t\\ & + 
   130092164100 s^6 t - 141668392500 s^7 t - 2988082710 t^2 + 
   40669277040 s t^2 - 143003556000 s^2 t^2\\ & + 138926113200 s^3 t^2 - 
   205071567345 s^4 t^2 + 34649953200 s^5 t^2 + 191909250 s^6 t^2 - 
   2983517550 t^3\\ & + 11422409400 s t^3 + 11281245840 s^2 t^3 + 
   12503166480 s^3 t^3 - 179607375 s^4 t^3\\ & - 456350355 t^4 - 
   1281756960 s t^4 + 54325080 s^2 t^4 - 334611 t^5 + 1180980 s t^5\Big)\end{aligned}&&\end{flalign}

\begin{flalign}\label{icodtfactor}\begin{aligned} &H_2(s,t)=\\&9 \Big(41065603 - 1295838420 s + 15320351025 s^2 - 84495071970 s^3 + 
   226929135060 s^4 - 311796618894 s^5\\ & + 180564347880 s^6 + 
   192831289200 s^7 - 2827791000 s^8 + 351078255 t - 7191203760 s t\\ & + 
   46956119580 s^2 t - 109054125630 s^3 t + 115785174600 s^4 t - 
   146492827020 s^5 t - 7724921400 s^6 t\\ & + 39366000 s^7 t + 
   864900990 t^2 - 8789339160 s t^2 + 13789435950 s^2 t^2 + 
   7464559050 s^3 t^2\\ & + 14634835380 s^4 t^2 - 75779550 s^5 t^2 + 
   640900350 t^3 - 865177200 s t^3 - 2466017460 s^2 t^3\\ & - 
   82471770 s^3 t^3 + 46987695 t^4 + 32411340 s t^4 + 
   688905 s^2 t^4 + 649539 t^5\Big)\end{aligned}&&\end{flalign}

\begin{flalign}\label{icointeqpoly1}\begin{aligned}&q_3(s)=\\&+4519637361671274496000 - 209181231007564352716800 s + 
  3321457506854369734164480 s^2\\& - 16912751468864437544288256 s^3 - 
  116921878612255069968531456 s^4\\& + 
  2027045152728857868094144512 s^5 - 
  10194970827460336472888967168 s^6\\& + 
  6976444965746627405502480384 s^7 + 
  172944934506481884571889565696 s^8\\& - 
  988250331686343727660081348608 s^9 + 
  2553828785423744289114090110976 s^{10}\\& - 
  2344640376723247467482041810944 s^{11} - 
  4966685189992395473902257881088 s^{12}\\& + 
  17576539179306720974557030416384 s^{13} - 
  18243873711150452680764558670848 s^{14}\\& - 
  31015621379102614675220323584 s^{15} + 
  5746571146290950737285450588416 s^{16}\\& + 
  219198145374369072919385534592 s^{17} - 
  39301095618302669628042931344 s^{18}\\& + 
  17859592862720046236039186304 s^{19} + 
  33006943093923122542733952 s^{20}\\& - 3158681972301971439233508 s^{21} + 
  13384995721483031782209 s^{22}\\& + 4408387278742778724 s^{23}\end{aligned}&&\end{flalign}

\begin{flalign}\label{icointeqpoly2}\begin{aligned}&q_4(t)=\\&+957577417423367365285503857386946905759044128019\\&  - 
 38193396090594794919315882235619599293980769162302t \\& + 
 599658788799156830172697892714918339582558538429801 t^2\\&  - 
 4455642269862328433662237645448301923361567265814496 t^3\\& + 
 11423440036029514236194427466474444430806461119719089 t^4\\&  + 
 56988573353050813089403803321177863192246321364031458 t^5\\& - 
 544247180911451122803482546072465816049457363080376461 t^6\\&  + 
 1780162762115049663091681919054561595334044575621593284 t^7\\& - 
 2605798650218188168518084298397156587624905489626697522 t^8 \\& + 
 495935564405864785383853094256048256898832235213364484 t^9\\& + 
 2510805036529709033931299336293637425233364385624374602 t^{10}\\&  - 
 1017253515010728770365790510399170286004994183849832208 t^{11}\\& - 
 149990141942512043522527232858411534917647005509572750 t^{12} \\& + 
 434320482118794710851486345331817104508976469564122452 t^{13} \\&- 
 1207671778700273053147945480687788072877247166920682 t^{14}\\&  - 
 1296602674662488339863162013313421699886037336607080 t^{15}\\& - 
 16215073118307783073925701402671068303411158821489 t^{16}\\&  + 
 108374307564297427167551141629985815480631124474 t^{17}\\& + 
 1413867303004562320958685439286796959982032925 t^{18}\\&  + 
 4948178294653340987297871157835191425144912 t^{19}\\& + 
 240377463274564672072528274760175543917 t^{20} \\& - 
 54117437609186157810550382898752502 t^{21} \\&- 
 1113371547932844477170246789721 t^{22}\\&  + 506404671871019220362052 t^{23}\end{aligned}&&\end{flalign}

\begin{flalign}\label{icoresultants}\begin{aligned}&\mathrm{Res}_s(G_2,H_2)(t)=\\&-55957321335751478724210392682094571010933699549277917763878154746093750000000000000000000000000000 \\ & \cdot(-361 + 3 t) (-73 + 3 t) ((-1 + 
   3 t)^2) (13 + 6 t) (53041 - 115926 t + 9 t^2)\\&\cdot (3521 - 2166 t + 
   729 t^2) (248369281 - 1475193612 t - 7509490506 t^2 - 
   238049388 t^3 + 81 t^4)\\&\cdot(3361 - 174252 t + 1784214 t^2 + 
   1101492 t^3 + 81 t^4)\\&\cdot (2545201 - 76777212 t + 578566854 t^2 + 
   6317892 t^3 + 81 t^4) \\&\cdot(423287 - 3795684 t + 11053773 t^2 + 
   68634 t^3 + 162  t^4)\, q_4(t) \end{aligned}&&\end{flalign}

\begin{flalign}\label{icoresultantt}\begin{aligned}&\mathrm{Res}_t(G_2,H_2)(s)=\\&18398940223618143886412514410156250000000000000000 (-4 + 3 s)(8 + 
   3 s) (1 + 6 s) (304 - 1212 s + 9 s^2)\\&\cdot (8 - 12 s + 45 s^2) (16 - 
   48 s + 81 s^2) (-2624 + 36384 s - 122904 s^2 - 16956 s^3 + 
   81 s^4)\\&\cdot (256 - 2496 s + 3456 s^2 + 7344 s^3 + 81 s^4) (256 - 
   11136 s + 16416 s^2 + 308664 s^3 + 81 s^4)\\&\cdot (784 - 8904 s + 
   27189 s^2 - 3024 s^3 + 324  s^4) \, q_3(s) \end{aligned}&&\end{flalign}

}

\setcounter{biburllcpenalty}{7000} 
\setcounter{biburlucpenalty}{8000} 
\setcounter{biburlnumpenalty}{9000}
\printbibliography[heading=bibintoc, title={References}]

\end{document}